\documentclass[a4paper,oneside,12pt]{article}

\usepackage{amsfonts}
\usepackage{amssymb}
\usepackage{amsmath}
\usepackage{latexsym}

\def\call#1{\ensuremath{\cal #1}}
\def\M#1{\ensuremath{\mathbb #1}}
\def\mfk#1{\ensuremath{\mathfrak #1}}
\def\f{\ensuremath{\varphi}}
\def\e{\ensuremath{\varepsilon}}
\def\wt#1{\ensuremath{\widetilde{#1}}}
\def\wh#1{\ensuremath{\widehat{#1}}}
\def\wb#1{\ensuremath{\overline{#1}}}
\def\z{{\bf z}}
\def\vol{{\rm vol}}
\def\svol{{\rm Strvol}}

\newenvironment{remark}{\begin{rem}\begin{rm}}{\end{rm}\end{rem}}

\newenvironment{defi}{\begin{defiteo}\begin{rm}}{\end{rm}\end{defiteo}}
\newenvironment{defitit}[1]{\begin{defiteo}[#1]\begin{rm}}{\end{rm}\end{defiteo}}

\newtheorem{teo}{Theorem}[section]
\newtheorem{prop}[teo]{Proposition}
\newtheorem{cor}[teo]{Corollary}
\newtheorem{lemma}[teo]{Lemma}
\newtheorem{defiteo}[teo]{Definition}
\newtheorem{rem}[teo]{Remark}
\newcommand{\proof}{\noindent{\em Proof.}\ }

\newcommand{\qed}{\\\makebox[\textwidth][r]{$\Box$}\vspace{\baselineskip}}

\begin{document}

\begin{center}
{\Large Hyperbolic volume of representations of fundamental groups of cusped
  $3$-manifolds \\\ \\ }

\vskip5mm
\textsc{Stefano Francaviglia\footnote{Author address: Stefano
      Francaviglia, Dipartimento di Matematica Applicata ``U. Dini'', 
      via Bonanno 25b,
      I-56126 Pisa, e-mail: s.francaviglia@sns.it.}\footnote{This work
    was partially written while the author was a PhD student at the
    Scuola Normale Superiore di Pisa, and was
    partially supported by the C.R.M. (Universit\'at Autonoma de Barcelona).}}

\today
\end{center}

\vskip5mm

\begin{abstract}
Let $W$ be a compact manifold and let $\rho$ be a representation of its
fundamental group into {\rm PSL}$(2,\M C)$. Then the volume of $\rho$ is
defined by taking any $\rho$-equivariant map from the universal cover
\wt W to $\M H^3$ and then by integrating the pull-back of the
hyperbolic volume form on a fundamental domain. It turns out that such
a volume does not depend on the choice of the equivariant map.
Dunfield extended this construction to the case of a 
non-compact (cusped) manifold $M$,
but he did not prove the volume is well-defined in all cases.

We prove here that the volume of a representation is always
well-defined and depends only on the representation. 
Moreover, we show that this volume can be easily computed by
straightening any ideal triangulation of $M$.

We show that the volume of a representation is bounded from above by the
relative simplicial volume of $M$.
Finally, we prove a rigidity theorem for representations of 
the fundamental group
of a hyperbolic manifold. Namely, we prove that if $M$ is hyperbolic and 
$\vol(\rho)=\vol(M)$ then $\rho$ is
discrete and faithful.

\end{abstract}

\section{Introduction}

Let $W$ be a compact manifold and let $\rho$ be a representation of its
fundamental group into {\rm PSL}$(2,\M C)\simeq$ Isom$^+(\M H^3)$. 
The volume of $\rho$ is
defined by taking any $\rho$-equivariant map from the universal cover
\wt W to $\M H^3$ and then by integrating the pull-back of the
hyperbolic volume form on a fundamental domain. 
This volume does not depend on the choice of the equivariant map
because two equivariant maps are always equivariantly homotopic and the
cohomology-class of the pull-back of the volume form is invariant
under homotopy. 

In \cite{D} this definition is extended to the case of a non compact cusped
$3$-manifold $M$ (see Definitions~\ref{d vdr} and~\ref{d pd}). When $M$
is not compact, some problems of integrability arise
if one tries to use the above definition of the volume of a
representation. 
The idea of Dunfield for overcoming these difficulties 
is to use a particular (and natural) class of equivariant maps, called
{\em pseudo-developing maps} (see Definition~\ref{d pd}), that have a
nice behavior on the cusps of $M$ allowing to control
their volume. Concerning the well-definition of the volume,
working with non-compact manifolds, 
two pseudo-developing maps in general are not equivariantly homotopic and
in~\cite{D} it is not proved that the volume of a representation does
not depend on the chosen pseudo-developing map.

\ 

In this paper we show that the volume of a representation is
well-defined even in the non-compact case, and 
we generalize to non-compact manifolds
some results know in the compact case. 
We restrict to the orientable case. 
The paper is structured as follows.

In Sections~\ref{s uno} and \ref{s due} we introduce the notion of 
pseudo-developing map for a given
representation $\rho:\pi_1(M)\to$ Isom$^+(\M H^3)$ and the
notion of straightening of such a map. 

In Section~\ref{s 3} 
we prove that for each orientable cusped $3$-manifold $M$ and 
for each representation $\rho:\pi_1(M)\to$ Isom$^+(\M H^3)$, the volume
of $\rho$ is well-defined and depends only on $\rho$. The main
theorems are:

\ \\
{\bf Theorem \ref{t vndpd}} {\em Let  $D_\rho$ and $F_\rho$ be two
  pseudo-developing maps for $\rho$. Then $\vol(D_\rho)=\vol(F_\rho)$.}

\ \\
{\bf Theorem \ref{t vsd}} {\em For any pseudo-developing map $D_\rho$ for
  $\rho$ we have $\vol(D_\rho)=$ $\svol(D_\rho)$.
}

\ 

Roughly speaking, Theorem~\ref{t vsd}  says that the volume of $\rho$
can be computed by straightening any ideal triangulation of $M$ and
then summing the volume of the straight version of the tetrahedra.

In Section~\ref{s qua}, generalizing the techniques used for the proof
of Theorem~\ref{t vsd}, we show that the volume of a
representation $\rho$ is bounded from above by the relative simplicial volume:

\ \\
{\bf Theorem \ref{t sv}} {\em For all representations
  $\rho:\pi_1(M)\to$ {\rm Isom}$^+(\M H^3)$ we have
  $|\vol(\rho)|<v_3\cdot||(\overline M,\partial \overline M)||$,
where $v_3$ is the volume of a regular ideal tetrahedron in $\M H^3$.}

\
 
In Section~\ref{s r}  we prove the following rigidity theorem for
representations of the fundamental group of a hyperbolic manifold:

\ \\
{\bf Theorem \ref{t ug}}
{\em   Let M be a non-compact, complete, orientable hyperbolic
  $3$-manifold of finite volume. Let $\Gamma\cong\pi_1(M)$ be the
  sub-group of ${\rm PSL}(2,\M C)$ 
  such that $M=\M H^3/\Gamma$. Let 
  $\rho:\Gamma\to$ {\rm PSL}$(2,\M C)$ be a representation.
If $|\vol(\rho)|=\vol(M)$ then $\rho$ is discrete and
faithful. More precisely there exists $\f\in{\rm PSL}(2,\M C)$ such
that for any $\gamma\in\Gamma$ 
$$\rho(\gamma)=\f\circ\gamma\circ\f^{-1}.$$
} 

In Section~\ref{s cor} we give some corollaries. 
In particular we show how from
Theorem~\ref{t ug} one can get a proof of Mostow's rigidity for non-compact
manifolds (see~\cite{P} and~\cite{BCS} for a more
general statement and a different proof):

\ \\
{\bf Theorem \ref{t mstw} (Mostow's rigidity for non-compact manifold)}
{\em Let $f:M\to N$ be a proper map between two orientable non-compact,
complete hyperbolic $3$-manifolds of finite volume. Suppose that
\vol(M)={\rm deg}(f)\vol(N).
Then $f$ is properly homotopic to a locally isometric covering with the
same degree as $f$.}

Other corollaries that can be
useful for checking the hyperbolicity of a $3$-manifold 
are also shown.

\vskip2ex
I would like to thank Carlo Petronio and Joan Porti for the very
many very interesting conversations about this theme.

\section{General definitions}\label{s uno}
We fix here the class of manifolds we consider, namely the class of ideally
triangulated cusped manifolds.
Since we work with
cusped manifolds, we want to fix a structure on the cusps.

\begin{defitit}{Cusped manifold}
  An orientable manifold $M$ is called {\em cusped manifold} if it is  
 diffeomorphic to the interior of a compact manifold with boundary
 $\wb M$. A {\em cusp} of $M$ is a closed regular neighborhood of a component
 of $\partial
 \overline M$. In the following we require $M$ to have dimension $3$ and
 $\partial \overline M$ to be a union 
 of tori, so each cusp is homeomorphic to $T^2\times[0,\infty)$.

We define $\wh M$ as the
  compactification of $M$ obtained by adding one point for each cusp of
  $M$.
  Called $\wt M$ the universal cover of $M$, we call $\wh{\wt M}$ the
  space obtained by adding to
  $\wt M$ one point for each lift of each cusp of $M$.

  We call the points added to $M$ (or $\wt M$) {\em ideal points} of $M$
  (or \wt M). For each ideal point $p$ of $M$, we fix a smooth product
  structure 
  $T_p\times[0,\infty)$ on the cusp relative to $p$. Such a structure
  induces a cone structure,  
  obtained from $T_p\times[0,\infty]$ by collapsing
  $T_p\times\{\infty\}$ to $p$, on a neighborhood $C_p$ of $p$ in $\wh M$. 

  We lift such structures to the universal cover. Let \wt p be an ideal
  point of \wt M that projects to the ideal point $p$ of $M$. 
  We denote by $N_{\wt p}$  the cone at \wt p. The cone
  $N_{\wt p}$ is homeomorphic to $P_{\wt p}\times[0,\infty]$ where
  $P_{\wt p}$ covers
  the torus $T_p$ and $P_{\wt p}\times\{\infty\}$ is collapsed to \wt p. 
\end{defitit}

\begin{remark}
In the definition of cusped manifold we have included a
fixed product structure on the cusps. This is for technical reasons,
however we will show that the results about the volume of representations
do not depend on the chosen structure.  
\end{remark}

\begin{remark}
  Let \wt M be the universal cover of $M$. In the following, when we
  speak about $\pi_1(M)$, we tacitly assume that a base-point and one
  of its lifts have been fixed. If $p$ is an ideal point of $M$, then
  $\pi_1(T_p)$ is well-defined only up to conjugation. Called $\{\wt
  p_i\}$ the set of the lifts of $p$, there is a one-to-one
  correspondence between the stabilizers Stab$(\wt p_i)$ of $\wt p_i$
  in the group of deck transformations of $\wt M\to M$ and the conjugates of
  $\pi_1(T_p)$ in $\pi_1(M)$. Such a correspondence is uniquely
  determined once the base-points have been fixed.
\end{remark}

To avoid pathologies, since we are working with cusped manifolds, 
we need that the maps we use have a nice
behavior ``at infinity.'' Namely, we will often require 
that a map from a cusp to $\M H^3$ is a cone-map in the following sense.

\begin{defitit}{Cone-map}
Let $A$ be a set, $c\in\M R$ and $C$ be the cone obtained from
$A\times[c,\infty]$ by collapsing $A\times\{\infty\}$ to a point,
which we call $\infty$. A map $f:C\to\wb{\M H}^n$ is a {\em cone-map}
if:
\begin{itemize}
\item $f(C)\cap\partial \M H^n=\{f(\infty)\}$; 

\item $\forall a\in A$ the map $f_{|a\times[c,\infty]}$ is either the
  constant to $f(\infty)$ or the geodesic ray from $f(a,c)$ to
  $f(\infty)$, parametrized in such a way that the parameter
  $(t-c),\ t\in[c,\infty]$, is the arc-length.
\end{itemize}
\end{defitit}

We recall here the definition of pseudo-developing map for a
representation (see \cite{D}). 

\begin{defitit}{Pseudo-developing map}\label{d pd} 
  Let $M$ be a cusped manifold and let  
 $\rho:\pi_1(M)\to{\rm Isom}^+(\M H^3)$ be a
  representation. 
A {\em pseudo-developing map} for $\rho$ is a piecewise smooth map
  $D_\rho:\wt M\to \M H^3$ which is equivariant w.r.t. the actions of
  $\pi_1(M)$ on \wt M via deck transformations and on $\M H^3$ via
  $\rho$. Moreover we require $D_\rho$ to extend to a continuous map, which we
  still call $D_\rho$, from $\wh{\wt M}$ to $\overline{\M H}^3$ that
  maps the ideal points to $\partial \M H^3$ (see Remark~\ref{r die} 
  for comments on this property). Finally we require that
  there exists $t_{D_\rho}\in \M R^+$ such that for each cusp
  $N_p=P_p\times[0,\infty]$ of \wt M, 
  the restriction of $D_\rho$ to $P_p\times[t_{D_\rho},\infty]$ is a cone-map.
\end{defitit}

Let $\gamma\neq{\rm id}$ be an isometry of $\M H^3$ and let Fix$(\gamma)$
be the set of fixed points of $\gamma$. Then
Fix$(\gamma)\cap\partial\M H^3$ consists of either one or two
points. Moreover, if $\gamma_1$ and $\gamma_2$ commute, then
Fix$(\gamma_1)$ is $\gamma_2$-invariant. It follows that if $\Gamma$
is an Abelian subgroup of orientation-preserving isometries and
$\gamma\in\Gamma$, 
then Fix$(\gamma)$ is $\Gamma$-invariant. 
Actually, for almost all Abelian $\Gamma$ and for any $\gamma_1,
\gamma_2\in\Gamma\setminus\{{\rm id}\}$ we have
$$\textrm{Fix}(\gamma_1)= \textrm{Fix}(\gamma_2).$$

The only cases in which this is not true are when $\Gamma$ is a
dihedral group generated by two rotations of angle $\pi$ around
orthogonal axes. Such a group is isomorphic to $\M Z_2\times \M Z_2$ and
its unique fixed point is the intersection of the axes. 
It follows that any Abelian group $\Gamma$ of orientation-preserving 
isometries has a fixed
point in $\wb{\M H}^3$ and, if $\Gamma$ is not dihedral, then it has a
fixed point in $\partial \M H^3$.

Let now $p$ be an ideal point of $\wh{\wt M}$. Since Stab$(p)$ is
Abelian, then either it is dihedral or it has a fixed point in
$\partial\M H^3$. If $\rho$ is a representation of $\pi_1(M)$ and
$D_\rho$ is a pseudo-developing map for $\rho$, then $D_\rho(p)$ is a
fixed point of $\rho($Stab$(p))$. 
It follows that, using Definition~\ref{d pd}, in
order for a pseudo-developing map to exist, $\rho($Stab$(p))$ must
have a fixed point in $\partial \M H^3$.

\begin{remark}~\label{r die}
We included in Definition~\ref{d pd} the requirement that $D_\rho$
maps ideal points to $\partial \M H^3$ only for
simplicity. No pathologies do occur if some ideal point is mapped to
the interior of $\M H^3$.
Coherently with this fact, from now on we suppose that:
\begin{itemize}
\item[] For each boundary torus $T$, the group $\rho(\pi_1(T))$ is not
  dihedral. 
\end{itemize}
As above we notice that this is only for simplicity and one can easily
check that all the results of this paper remain true, {\em mutatis
mutandis},  without this assumption.
\end{remark}

\begin{lemma}\label{l 1}
  Let M be a cusped manifold and let
  $\rho:\pi_1(M)\to$ {\rm Isom}$^+(\M H^3)$ be a representation. Then a
  pseudo-developing map $D_\rho$ exists.
\end{lemma}

\proof The proof is the same as in \cite{D}, we recall it by completeness.
We construct a pseudo-developing map inductively on the
$n$-skeleta. Let $p$ be an ideal point of \wt M. Since Stab$(p)$ is
Abelian and not dihedral, then its $\rho$-image has at least one fixed point
$q\in\partial \M H^3$. We define $D_\rho(p)=q$ and, for all
$\alpha\in\pi_1(M)$ we set $D_\rho(\alpha(p))=\rho(\alpha)(q)$.
We do the same for the other ideal points.
Now for each ideal point $p$ we define $D_\rho$ on $P_p\times\{0\}$ in
any Stab$(p)$-equivariant way and then we make the cone over
$D_\rho(p)$ in such a way that $D_\rho$ has the cone property.
Then we extend $D_\rho$ in any
equivariant way. The extension is possible because $\M H^3$ is contractible.
\qed
\begin{remark}\label{r pip} Let $p$ be an ideal point of $\wh{\wt M}$.
  If $\rho($Stab$(p))$ is a parabolic non-trivial group, then it has a
  unique fixed point. It follows that $D_\rho(p)$ is uniquely
  determined.   
  Thus, if all the $\rho$-images of the stabilizers of the ideal
  points are parabolic, then the $D_\rho$-images of all the ideal points
  are uniquely determined.  
\end{remark}

\begin{defitit}{Ideally triangulated manifold}
 Let $M$ be a cusped manifold.
 An ideal triangulation of $M$ is a triangulation of $\wh M$ having
 the set of ideal points as $0$-skeleton. 
 An {\em ideally triangulated manifold} is a cusped manifold equipped
 with a finite smooth ideal triangulation $\tau$.
 We require the triangulation to be compatible with the product
 structure. 
That is, for each cusp $N_p$ we require $\tau\cap (T_p\times\{0\})$ 
to be a triangulation of $T_p$ and the restriction
  to $N_p$ of  $\tau$ to be the product triangulation.
\end{defitit}

We will often consider the simplices of an ideal triangulation of a manifold
$M$ as subsets of $\wh M$.
\begin{remark}\label{r acmcbid} It is well-known that any cusped manifold 
 can be ideally triangulated (see for example~\cite{BP}). 
\end{remark}

\section{The straightening}\label{s due}

A straightening of a pseudo-developing map $D$ is a map that agrees
with $D$ on the ideal points and that maps each
tetrahedron to a straight one. The straightening is useful to
calculate the hyperbolic volume associated to a pseudo-developing map
(see Section~\ref{s 3}).
A particular case is when the manifold
$M$ is complete hyperbolic, because in this case the straightening descends to
a map from $M$ to itself. Here we prove that such a map is onto. 

Let $\Delta\subset\M H^3$ be an oriented geodesic ideal tetrahedron. Since
$\Delta$ is the convex hull of its vertices, then the Isom$^+(\M
H^3)$-class of $\Delta$ is completely determined by the Isom$^+(\M
H^3)$-class of the oriented set of its vertices (the orientation of
the vertices is defined up to the action of ${\call A}_4$). Such a
class is completely 
determined by a non-real complex 
number called modulus, up to a three-to-one
ambiguity. Such an ambiguity
can be avoided by choosing a preferred pair of opposite edges of
$\Delta$ (see \cite{BP}, \cite{F3}, \cite{PP}, \cite{PW} \cite{Th1}). 
We extend the notion of modulus to
the set of flat 
tetrahedra, that is to those whose vertices are distinct and lie on a
hyperbolic plane of 
$\M H^3$, by accepting real moduli different from 0 and 1.  We
want to extend this definition also to the degenerate tetrahedra,
{\em i.e.} to those having two ore more coincident vertices. Unfortunately, for
such a tetrahedron it is not possible to encode its isometry class in a
complex number. Let us agree that when we use a modulus in
$\{0,1,\infty\}$ for $\Delta$, we mean that $\Delta$ is a degenerate
tetrahedron and that the modulus encodes the complete information on
the isometry class of $\Delta$, {\em i.e.} who are the coincident
vertices of $\Delta$. 

\begin{defi} Let $\Delta^k$ be the standard $k$-simplex.
Let $\f:\Delta^k\to \wb{\M H}^n$ be a continuous map that maps
the $0$-skeleton of $\Delta^k$ to $\partial \M H^n$. 
Let $Q$ be the
Euclidean convex hull of the \f-image of the vertices  of
$\Delta^k$, made in a disc model of $\M H^n$. 
Let $\psi:\Delta^k\to Q$ be the only
simplicial map that agrees with $\f$ on the
$0$-skeleton.

We say that the map $\f$ is {\em standard} if there exist two
homeomorphisms $\eta:{\rm Im}(\f)\to Q$ and $\beta:\Delta^k\to \Delta^k$
such that $$\eta\circ\f\circ\beta=\psi.$$
We say that a foliation $\cal F$ of $\Delta^k$  is {\em standard} if
there exists a standard map $\f:\Delta^k\to \wb{\M H}^n$ such that ${\call
F}=\{\f^{-1}(x)\}$.  
\end{defi}

\begin{remark}\label{r uno}
  For any standard map $\f$ the dimension of ${\call F}=\{\f^{-1}(x)\}$ 
  depends only on the \f-image of the $0$-skeleton.
\end{remark}

\begin{remark}
 It is not hard to show that for a map  \f\ 
 to be standard does not depend on the disc model we use. In other
 words \f\ is standard if and only if $\gamma\f$ is standard for any
 isometry $\gamma$.
\end{remark}

Let $M$ be an ideally triangulated manifold,
$\rho:\pi_1(M)\to${\rm Isom}$^+(\M H^3)$ be a representation, and $D_\rho$ be
a pseudo-developing map for $\rho$. Let $\Delta$ be a tetrahedron of
$\tau$ and \wt\Delta\ be one of its lifts. The vertices of \wt\Delta\
are ideal points, so their images under $D_\rho$ lie in $\partial
\M H^3$. The map $D_\rho$ determines a modulus for \wt\Delta{ } simply
by considering the convex hull of the image of its vertices (the
orientation is the one induced by $M$). Note that
since $D_\rho$ is equivariant, then it defines a modulus for
$\Delta$. For each face $\sigma$ of $\Delta$ we call
${\rm Str}_{D_\rho}(\wt\sigma)$, or simply ${\rm Str}(\wt\sigma)$, 
the straight simplex obtained as the convex hull of the $D_\rho$-image
of the vertices of \wt\sigma.

\begin{defitit}{Straightening}\label{d str}
  A straightening of $D_\rho$ is a continuous, piecewise smooth, 
 $\rho$-equivariant map
  ${\rm Str}(D_\rho):\wh{\wt M}\to 
  \overline{\M H}^3$ such that:
  \begin{enumerate}
   \item For each simplex $\sigma$ of the triangulation,
     ${\rm Str}(D_\rho)$ maps \wt\sigma{ } to 
     ${\rm Str}(\wt\sigma)$. 
   \item The restriction of ${\rm Str}(D_\rho)$ to any simplex
   $\sigma$ is standard.
   \item For each cusp $N_{\wt p}=P_{\wt p}\times[0,\infty]$ there
     exists $c\in\M R$ such that ${\rm Str}(D_\rho)$ restricted to
     $P_{\wt p}\times[c,\infty]$ is a cone-map.
\end{enumerate}
\end{defitit}

\begin{lemma}\label{l 2}
  Let M be an ideally triangulated manifold. Let $\rho$ be a representation
  $\rho:\pi_1(M)\to${\rm Isom}$^+(\M H^3)$ and $D_\rho$
  be a pseudo-developing map. Then a straightening ${\rm Str}(D_\rho)$ of
  $D_\rho$ exists. Moreover ${\rm Str}(D_\rho)$ is always equivariantly
  homotopic to $D_\rho$ via a homotopy that fixes the ideal points.
\end{lemma}
\proof
A straightening of $D_\rho$ can be constructed with the same techniques of
Lemma~\ref{l 1}. Regarding the homotopy, since $D_\rho$
maps non-ideal points to the interior of $\M H^3$, then one can use a
geodesic flow  with the time-parameter in $[0,\infty]$ 
(for example the convex combination of
Definition~\ref{d cvxc}) to construct a
homotopy with the required properties.\qed
\begin{remark}
A straightening in general is not a pseudo-developing map in our
setting, because it can map some point of \wt M to $\partial \M
H^3$. However, if there are no degenerate tetrahedra, then 
a straightening is also a pseudo-developing map, and the
homotopy between $D_\rho$ and ${\rm Str}(D_\rho)$ can be made coherently
with the cone structure of the cusps, {\em i.e.} in such a way that the
intermediate maps along the homotopy 
between $D_\rho$ and ${\rm Str}(D_\rho)$ have the cone property on the cusps.
\end{remark}

When $M$ has a complete hyperbolic structure of finite volume,
there is a natural notion of straightening of the ideal
triangulation.
Namely, choose the arc-length as the cone
parameter on the cusps of $M$ and consider $\M H^3$
as the universal cover of $M$.
 Then choose $\rho$ as the holonomy of the
hyperbolic structure of $M$;  the identity map of $\M H^3$ clearly is
a pseudo-developing map for $\rho$. A natural straightening map is a
straightening of the identity.  

\begin{prop}\label{p so}
  Let M be an ideally triangulated manifold equipped with a complete,
  finite-volume hyperbolic structure. Then any natural straightening
  map projects to a map ${\rm Str}:\wh M\to \wh M$ which is onto. Moreover
  ${\rm Str}(M)\supset M$.
\end{prop}

\proof
It is easy to see that $\wh{\wt M}=\wh{\M H^3}$ naturally embeds into
$\overline{\M H}^3$ and that the ideal points lie on $\partial \M
H^3$. Since the straightening is equivariant, then it projects to a
map ${\rm Str}:\wh M\to \wh M$. Moreover, Str fixes the ideal points. We
prove that Str is onto. One can easily prove that $H_3(\wh M;\M
Z)\cong H_3(\overline M,\partial \overline M;\M Z)\cong\M Z$. So we can define
the degree of a map $f:\wh M\to\wh M$ by
$$f_*([\wh M])=\textrm{deg}(f)\cdot[\wh M]$$
where $[\wh M]$ is the generator of $H_3(\wh M;\M Z)$ induced by the
orientation of $M$.
Now note that by Lemma~\ref{l 2} the natural straightening is
homotopic to the identity via an equivariant homotopy. Because of
equivariance, the homotopy projects to a homotopy between Str and the
identity. It follows that ${\rm Str}_*$ and id$_*$ coincide on $H_*(\wh M;\M
Z)$, so $\textrm{deg}(\textrm{Str})=\textrm{deg}(\textrm{id})=1$.
Now suppose that Str is not onto and let $x$ be a point in $\wh M$
outside its image. If we consider Str as a map from \wh M to $\wh
M\setminus \{x\}$, we get ${\rm Str}_*([\wh M])=0\in H_3(\wh M\setminus
\{x\};\M Z)$ simply because $H_3(\wh M\setminus\{x\};\M Z)=0$. Then
${\rm Str}_*([\wh M])$ is a boundary in $\wh M\setminus \{x\}$, consequently it is a
boundary also in \wh M. It follows that ${\rm Str}_*([\wh
M])=0$. This implies $\textrm{deg}(\textrm{Str})=0$, that is a
contradiction. 

The last assertion follows because Str is onto and fixes
the ideal points. 
\qed

\section{Volume of representations}\label{s 3}

For this section we fix an ideally triangulated manifold
$M$ and a representation $\rho:\pi_1(M)\to{\rm Isom}^+(\M H^3)$.

In this section we recall the notion of volume of an equivariant map from \wt M
to $\M H^3$. 
We prove that if we restrict to the class of pseudo-developing
maps, then the volume of $\rho$ is well-defined. 
Namely the volume does not
depend neither on the pseudo-developing map nor on the product
structure of the cusps.
Such a volume can be calculated
using a straightening of any 
pseudo-developing map and it is exactly the algebraic sum of
the volumes of the straightened tetrahedra.

\begin{defitit}{Volume of pseudo-developing map}\label{d vdr}
  Let $D_\rho$
  be a  pseu\-do-de\-vel\-op\-ing map for $\rho$. Let $\omega$ be the volume
  form of $\M H^3$ and let $D_\rho^*\omega$ be the pull-back of
  $\omega$. Since $D_\rho $ is equivariant, then $D_\rho^*\omega$
  projects to a $3$-form, that we still call $D_\rho^*\omega$, on
  $M$. The volume $\vol(D_\rho)$ of $D_\rho$ is defined by:
$$\vol(D_\rho)=\int_M D_\rho^*\omega.$$ 
\end{defitit}

\begin{remark}
We will see
below that for pseudo-developing maps the volume is always finite.
The same definition of volume does not work for any equivariant map
from \wt M 
to $\M H^3$ because if the pull-back of the
volume 
form is not in $L^1$, then the volume is not well-defined. 
\end{remark}

\begin{defitit}{Straight volume}\label{d strvol}
  Let $D_\rho$ be a pseudo developing map for $\rho$. Let $\{\Delta_i\}$
  be the set of the tetrahedra of the ideal triangulation of $M$ and
  $\{\wt\Delta_i\}$ be a set of lifts of the $\Delta_i's$. Let
  $v_i=0$ if ${\rm Str}(\wt\Delta_i)$ is a degenerate tetrahedron, and 
  let $v_i$ be the algebraic volume of ${\rm Str}(\wt\Delta_i)$ otherwise. 

  We define the straight volume of $(D_\rho)$ as 
  $\svol(D_\rho)=\sum_i v_i$.
\end{defitit}

\begin{remark} Let
$\Delta_1,\dots,\Delta_n$ be the tetrahedra of $\tau$. Let
$\z=(z_1,\dots,z_n)\in\{\M C\setminus\{0,1\}\}^n$ be a solution of
Thurston's hyperbolicity equations (see~\cite{BP},
\cite{F3}, \cite{NZ}, \cite{PP}, \cite{PW}, \cite{Th1}). Then there exists
a developing map $D_\z:\wt M\to \M H^3$ for \z\ that is a
pseudo-developing map for some holonomy $\rho(\z)$. Such a map is
already straight and we have
$\svol(D_\z)=\vol(D_\z)=\sum v_i=\vol(\z)$, where $v_i$ is
the volume of 
the geodesic ideal tetrahedron of modulus $z_i$.      
\end{remark}

Let $C_p=T_p\times[0,\infty]/_\sim $ be a cusp of $M$ and let
$N_{\wt p}=P_{\wt p}\times [0,\infty]/_\sim$ be one of its lifts in \wt M. 
Then we can identify $\pi_1(C_p)$ with ${\rm Stab}(\wt p)$.
Let $f:P_{\wt p}\times\{0\} \to \M H^3$ be a
${\rm Stab}(\wt p)$-equivariant map, let 
$\xi\in\partial\M H^3$ be a fixed point of
$\rho({\rm Stab}(\wt p))$ and let $F:N_{\wt p}:\to\overline{\M H}^3$ be
the cone-map obtained by coning $f$ to $\xi$.
As above, let
$F^*\omega$ be the pull-back of the volume-form on $C_p$. Similarly we
can pull-back the metric. We call $A^p_t$ the area of the torus
$T_p\times\{t\}$.  

\begin{lemma}\label{l at}
In  the previous setting, for $t>r$ we have:
$$ A^p_t\leq A^p_r e^{-(t-r)}\qquad \textrm{ and } \qquad \int_{T_p\times
  [t,\infty)} |F^*\omega|\leq A^p_t.$$
\end{lemma}

\proof
Let $(x,y)$ be local coordinates on $P_{\wt p}$.
Choose the half-space model $\M C\times \M R^+$ of $\M H^3$ and assume
that $\xi=\infty$. In such a model the hyperbolic metric at the
point $(z,s)$ is the
Euclidean one rescaled by the factor $1/s$. It follows that, called
$\alpha+i\beta$ and $h$ the complex and real components of $F$, we have
$$\alpha(x,y,t)+i\beta(x,y,t)=\alpha(x,y,r)+i\beta(x,y,r) 
\qquad h(x,y,t)=h(x,y,r)e^{(t-r)}.$$

The element of area at level $t$ is
  $d\sigma_t(x,y)=\sqrt{\det(^TJF_t\cdot H\cdot 
  JF_t)}$, where $F_t$ is the restriction of $F$ to $P_{\wt p}\times\{t\}$ and
  $H(x,y,t)=\frac{1}{h^2}{\rm Id}$ is the matrix 
  of the hyperbolic metric. From direct calculations it follows that
  $d\sigma_t(x,y)\leq d\sigma_r(x,y)e^{-t+r}$ and the first inequality
  follows. 

Now note that the volume element $|F^*\omega|$ at the point $(x,y,t)\in
C_p$ is bounded by the area element of the torus $T_p\times\{t\}$
multiplied by the length element of the ray $\{(x,y)\}\times[0,\infty]$. Since
the parameter $t$ is exactly the arc-length, then the length element is
exactly $dt$. It follows that
$$ \int_{T_p\times[t,\infty)}|F^*\omega|\leq\int_t^\infty
A^p_sds\leq\int_t^\infty A^p_t e^{-(s-t)}ds=A^p_t.$$
This completes the proof.
\qed
\begin{remark}
From Lemma~\ref{l at} it follows in particular that 
$\int_{T_p\times [t,\infty)} |F^*\omega|\leq A^p_0e^{-t}$. This
means that we have an estimate of $\int_{T_p\times [0,\infty)} |F^*\omega|$
not depending on the point $\xi=F(p)$ but only on the area of $T_p\times\{0\}$.
\end{remark}

\begin{remark} From Lemma~\ref{l at} it follows that $\vol(D_\rho)$ is
  finite for any pseu\-do-de\-vel\-op\-ing map $D_\rho$. 
\end{remark}

The following lemma is proved in \cite{D}.

\begin{lemma}\label{l f1} If $D_\rho$ and $F_\rho$ are two
  pseudo-developing maps for $\rho$ that agree on the ideal points,
  then $\vol(D_\rho)=\vol(F_\rho)$.  
\end{lemma}
This is because any two pseudo-developing maps are equivariantly
homotopic. The fact 
that they coincide on the ideal points allows one to construct a
homotopy $h$ that respects the cone structures of the cusps. Namely,
for each ideal point \wt p of \wt M we choose any equivariant homotopy
between the restrictions of $D_\rho$ and $F_\rho$ to $P_{\wt
  p}\times\{\bar t\}$, where $\bar t=$ max$\{t_{D_\rho},t_{F_\rho}\}$,
we cone such a
homotopy to $D_\rho(\wt p)$ along 
geodesic rays, and we extend the homotopy outside the cusps in any
equivariant way.    
For such a homotopy $h$ we can use the Stokes theorem on $M\times[0,1]$ for
$h^*\omega$ to obtain the thesis. More precisely, let 
$K_t$ be $M\setminus\cup_p
(T_p\times(t,\infty))$, where $p$ varies on the set of the ideal points; then
we have 
$$
\label{e 1}0=\int_{K_t\times[0,1]}d(h^*\omega)=
\int_{\partial(K_t\times[0,1])}h^*\omega= \int_{K_t}
(D_\rho^*\omega-F_\rho^*\omega) + \int_{\partial K_t\times[0,1]}
h^*\omega
$$
and, as in Lemma~\ref{l at}, we can prove that the last integral goes
to zero as $t\to \infty$.

We now prove that the claim of Lemma~\ref{l f1} is true in general.  

\begin{teo}\label{t vndpd}
Let  $D_\rho$ and $F_\rho$ be two
  pseudo-developing maps for $\rho$. Then $\vol(D_\rho)=\vol(F_\rho)$.  
 \end{teo}

\proof For $t\in[0,\infty)$, let $D_\rho^t$ be the map constructed as
follows: $D_\rho^t$ coincides with $D_\rho$ until the level $t$ of each
cusp. Then for each cusp $N_p$ we complete $D_\rho^t$ by coning 
 $D|_{P_p\times\{t\}}$ to $F_\rho(p)$ along geodesic rays in such a
 way that the arc-length is the  
 parameter $s-t$, where $s\in[t,\infty)$. Now, $D^t_\rho$ is a
 pseudo-developing map that agrees with $F_\rho$ on the ideal points. Thus by
 Lemma~\ref{l f1}    
$\vol(D^t_\rho)=\vol(F_\rho)$.
Since $D^t_\rho$ and $D_\rho$ agree outside the cusps and where they
differ they are cones on the same basis (and different vertices), 
from Lemma~\ref{l at} we get 
$$|\vol(D_\rho)-\vol(D_\rho^t)|\leq 2\sum_p
A_t^p\leq2(\sum_p A_0^p)e^{-t}$$
where $p$ varies on the set of ideal points and $A^p_t$ is the area
of the torus $T_p\times\{t\}$. As $t\to \infty$ we get the thesis.\qed

Similar techniques actually allow to prove the following theorem.
\begin{teo}\label{t vsd}
For any pseudo-developing map $D_\rho$ for $\rho$ we have $\vol(D_\rho)=$
 $\svol(D_\rho)$.
\end{teo}
Before proving Theorem~\ref{t vsd}, we give the following definition.

\begin{defi}\label{d cvxc}
  Let $f,g$ be two maps from a set $X$ respectively to $\M H^n$ and
  $\overline{\M H}^n$. For $t\in[0,\infty]$ the convex combination
  $\Phi_t$ from $f$ to $g$ is defined by: 
$$\Phi_t(x)=\left\{\begin{array}{lcl}
\gamma_x(t) &\ & t\leq{\rm dist}(f(x),g(x))\\
g(x)&\ &t\geq {\rm dist}(f(x),g(x))
\end{array}\right.$$  
where $\gamma_x$ is the geodesic from $f(x)$ and $g(x)$, parametrized
by arc-length.
\end{defi}
\begin{remark}
In Definition~\ref{d cvxc}, if $X$ is a topological space and $f$ and
$g$ are continuous, then the 
convex combination from $f$ to $g$ is continuous on
$X\times[0,\infty]$ because the function ${\rm dist}(f(x),g(x))$ is
well-defined and continuous from $X$ to $[0,\infty]$. 
\end{remark}

\noindent{\em Proof of \ref{t vsd}.}
For the proof assume that $t_{D_\rho}=0$. We start by fixing a suitable
homotopy $h$ between $D_\rho$ and ${\rm Str}(D_\rho)$. 
Define $h:\wt M\times [0,\infty]\to \M H^3$ outside the cusps to be the convex
combination from $D_\rho$ to ${\rm Str}(D_\rho)$ and then for each cusp 
$N_{\wt p}$ extend $h$ by coning 
$h((x,0),s)$ to $D_\rho(\wt p)$ along geodesic rays in such a way that
the parameter $t\in[0,\infty)$ of the cusp is the arc-length. 
Let $D_s(x)=h(x,s)$. By Lemma~\ref{l f1} we have that 
$$
\int_MD_\rho^*\omega=\int_MD_s^*\omega\qquad\textrm{ for } s\in(0,\infty).
$$
So we only have to prove that $\int_MD_s^*\omega\to
\svol(D_\rho)$ as $s\to \infty$. 
Clearly, it suffices to prove that for any tetrahedron $\Delta$ we have
$\int_\Delta D_s^*\omega\to v$ where $v$ is the volume of ${\rm Str}(\Delta)$.
If $\Delta$ does not collapse in the straightening, then the distance
from $D_\rho$ and ${\rm Str}(D_\rho)$ is bounded outside the cusps and so
$D_s={\rm Str}(D_\rho)$ for $s>>0$; since ${\rm Str}(D_\rho)$ is a
homeomorphism on $\Delta$, then $\int_\Delta D_s^*\omega$ is exactly
the volume of the straight version of $\Delta$.

If $\Delta$ collapses in the straightening, then we have to show that 
$\int_\Delta D_s^*\omega\to 0$.
This follows from direct calculations. We give only the lead-line of
them because they are involved but use elementary
techniques. Moreover, in the next section, we will give an alternative
proof of this theorem (see Theorem~\ref{t sv} and Remark~\ref{r difp}).

Given the convex combination $\Phi_t$ from a map $f$ to a map $g$, it
is possible to calculate the Jacobian of $\Phi_t$ as a function of
the derivatives of $f$ and $g$, the time $t$ and the distance between  
$f$ and $g$. This is not completely trivial, for
example think  of a tetrahedron as a convex combination of two segments: the
segments have zero area but in the middle we have quadrilaterals with
non-zero area. Using these calculations, we can estimate
$|D^*_s\omega|$ outside the cusps, showing that its integral goes to
zero as $s$ goes to infinity. Looking inside the cusps, by 
Lemma~\ref{l at} we reduce the estimate to the same estimate as above,
made with 2-dimensional objects (the bases of the cusps).\qed 

\begin{remark}
  Since $\vol(D_\rho)=\svol(D_\rho)$ it 
  follows that such a volume does not depend on the chosen cone
  structure of the cusps. 
  Moreover, by Theorem~\ref{t vndpd}, $\vol(D_\rho)$ 
 does not depend on the developing map, but only on $\rho$. This
  allows us to give the following definition.
\end{remark}

\begin{defi} The volume $\vol(\rho)$ of $\rho$ is the volume of any
  pseu\-do-de\-vel\-op\-ing map for $\rho$.
\end{defi}

As the following corollary shows, for hyperbolic manifolds the volume
of the holonomy is exactly the hyperbolic volume. 

\begin{cor}\label{c nat}
  Let $M$ be a complete hyperbolic manifold of finite volume. If
 $\rho$ is the holonomy of  the hyperbolic structure then $\vol(\rho)=$
 $\vol(M)$. 
\end{cor}
\proof
Consider $\M H^3$ as the universal cover of $M$ and choose the arc
length as the cone parameter of the cusps. 
Clearly the identity of
$\M H^3$ is a pseudo-developing map for $\rho$. 
Obviously we have $\int_M {\rm Id}^*(\omega)=$ $\vol(M)$.\qed

\begin{cor}
Let $z_i$ be the modulus induced by a pseudo-developing map $D_\rho$ on
the i$^{th}$ tetrahedron and let $v_i$ be the volume of a hyperbolic
ideal geodesic tetrahedron of modulus $z_i$. Then we have $\vol(\rho)=\sum v_i$.
\end{cor}

\begin{remark}
Even if $\sum v_i$ depends only on $\rho$, the moduli $z_i$ induced by
a pseudo-developing map $D_\rho$ actually can depend on  $D_\rho$. Namely,
any ideal point $p$ is mapped to a fixed point of $\rho($Stab$(p))$ and,
if this is not a parabolic group, we have more than one possibility for
$D_\rho(p)$. Conversely, if each $\rho($Stab$(p))$ is a non-trivial
parabolic group, then by Remark~\ref{r pip} it follows that the moduli
$z_i$ are uniquely determined by $\rho$.
\end{remark}

\begin{prop}\label{p pos}
  Let $g$ be a reflection of $\M H^3$ and let $\wb \rho$ be the
  representation $g\circ\rho\circ g^{-1}$. Then $\vol(\wb \rho)=-\vol(\rho)$.
\end{prop}
\proof If $D_\rho$ is a pseudo-developing map for $\rho$, then $g\circ
D_\rho$ is a pseudo-developing map for $\wb\rho$ and it is easily
checked that $\vol(g\circ D_\rho)=-\vol(D_\rho)$.\qed

We recall here two facts proved in \cite{D}.
\begin{prop}
  Suppose that $\rho_t:\pi_1(M)\to{\rm Isom}^+(\M H^3)$, 
$t\in[0,1]$ is a smooth one parameter family
  of representations. Then
  $\vol(\rho_0)=\vol(\rho_1)$.   
\end{prop}
This can be proved by using the one-parameter family  $\rho_t$ to
construct a one-parameter family of pseudo-developing maps. Some
estimates on the derivative along the $t$-direction of the volume of
such maps are needed.  

\begin{prop}\label{p f2}
  Suppose that $\rho$ factors through the fundamental group of a Dehn
  filling N of M. Then the volume of $\rho$ w.r.t. N
  coincides with the volume of $\rho$ w.r.t. M. 
\end{prop}

Theorem~\ref{t vsd} extends from ideal to ``classical''
triangulations, namely to genuine triangulations $\cal T$ of
$\overline M$. 
Consider such a \call T as a triangulation of $M$ with some simplices
at infinity (those in $\partial \overline M$).
Given a pseudo-developing map $D_\rho$ for $\rho$, define a
straightening of $D_\rho$ relative to \call T, exactly as in
Section~\ref{s due}, by considering the convex hulls of the
images of the vertices of \call T. 
Then, one can give 
the definition of the straight volume relative to \call T of a
developing map $D_\rho$ exactly as in Definition~\ref{d strvol}, with
the unique difference that one has to use the 
tetrahedra of \call T instead of
the ideal tetrahedra of an ideal triangulation of $M$. Call such a
volume $\svol^{\call T}(D_\rho)$.

Finally, exactly as in Theorem~\ref{t vsd}, one can prove the following
fact:
\begin{prop}\label{p vsd}
  Let \call T be a triangulation of $\overline M$ and $D_\rho$ be a
  pseudo-developing map for $\rho$. 
Then $\vol(\rho)=\svol^{\call T}(D_\rho)$.
\end{prop}

\section{Comparison with simplicial volume}\label{s qua}

Here we generalize the argument used to prove Theorem~\ref{t vsd} 
to compare $\vol(\rho)$ with the
simplicial volume of $M$, obtaining exactly the expected inequality.

Let $||(\overline M,\partial \overline M)||$ 
be the simplicial volume of $\overline M$ relative to the boundary
(see~\cite{BP}, \cite{G}, \cite{Ku}, \cite{Th1} for more 
details), 
and let $v_3$ be the volume of a regular straight ideal tetrahedron
of $\M H^3$.

\begin{teo}\label{t sv}
For any representation $\rho:\pi_1(M)\to$ {\rm Isom}$^+(\M H^3)$ we have
$$|\vol(\rho)|\leq v_3\cdot||(\overline M,\partial \overline M)||.$$ 
\end{teo}
\proof 
For the proof we fix 
a representation $\rho:\pi_1(M)\to{\rm Isom}^+(\M H^3)$ and
a pseudo-developing map $D_\rho$ for $\rho$.

Let $c=\sum_i\lambda_i\sigma_i$ be a smooth singular chain in 
$\overline M$; here each simplex $\sigma_i$ is a piecewise smooth map from 
the standard
tetrahedron $\Delta^3$ to $\overline M$. 
The simplicial volume of $c$ is defined as $||c||=\sum|\lambda_i|$. The
relative simplicial volume of 
$(\overline M,\partial \overline M)$ is defined as

$$
||(\overline M,\partial \overline M)||=\inf\{||c||:\, [c]=[\overline M]\in
H_3(\overline M,\partial \overline M)\}.
$$

The proof has two main steps:
\begin{enumerate}
\item\label{stp1} Given a smooth cycle
$c=\sum_i \lambda_i\sigma_i$ representing $[\overline M]$, 
show that $\vol(\rho)=\sum_i\int_{\Delta^3} \lambda_i
\sigma_i^*(D_\rho^*\omega)$, where $\omega$ is the volume form of $\M H^3$.
\item\label{stp3} By replacing $c$ with its straightening, show that
  $\vol(\rho)=\sum_i \lambda_iv_i$, where $v_i$ is the volume of a straight
  version of $\sigma_i$.
\end{enumerate}

From Step~\ref{stp3}  it follows that 
$$|\vol(\rho)|\leq\sum_i|\lambda_i|\cdot|v_i|\leq v_3\cdot||c||.$$
Theorem~\ref{t sv} follows taking to the infimum over all $c$'s
representing $[\overline M]$.

\vskip1ex\noindent{\em Step~\ref{stp1}}. Since a
pseudo-developing map has the cone property on the cusps, 
the $3$-form $D_\rho^*\omega$ 
defined on $M$ extends to a $3$-form on $\overline M$ that vanishes at
the boundary. So we can consider the class $[D_\rho^*\omega]\in
H^3(\overline M,\partial \overline M)$. Since $[c]=[\overline M]$, then 
$$
\vol(\rho)=\int_M D_\rho^*\omega=\langle [D_\rho^*\omega],
[\wb M]\rangle=\langle 
[D_\rho^*\omega],[c]\rangle=\sum_i\int_{\Delta^3}
\lambda_i\sigma_i^*(D_\rho^*\omega)
.$$ 

\vskip1ex\noindent{\em Step~\ref{stp3}}. 
The idea is the following. Consider a lift $\wt c$ of $c$ to $\wh{\wt
  M}$. Let $\wb c=(D_\rho)_*\wt c$ be the push-forward of $\wt c$ to 
$\M H^3$ via $D_\rho$ and let ${\rm Str}(\wb c)$ be a straightening of $\wb
c$. Since the straightening is homotopic to the identity, then there
exists a chain-homotopy of degree one {\em i.e.} a map $H$ from $k$-chains
to $(k+1)$-chains such that 
${\rm Str} - {\rm Id}=H\circ\partial - \partial\circ H$. Then we have
$$\begin{array}{lcl}
\vol(\rho)&=&\langle D_\rho^*\omega,\wt c\rangle=
\langle\omega, (D_\rho)_*\wt c\rangle= 
\langle\omega,\wb c\rangle\\&=&
 \langle\omega,\textrm{Str}(\wb c)\rangle +
\langle\omega,\partial H\wb c\rangle-
 \langle\omega,H\partial \wb c\rangle\\ &=&
\displaystyle{\sum_i}\lambda_iv_i+\langle d\omega,H\wb c\rangle- 
\langle\omega,H \partial \wb c\rangle.\end{array}$$ 
The last two summands are zero because $d\omega=0$ and, even if
$\partial \wb c\neq0$, everything can be made $\rho$-equivariantly 
so that the action of $\rho$ cancels out in pairs the contributions of
$\langle\omega,H \partial \wb c\rangle$. 

We formalize now this argument.
Let $C_k(X)$ denote the real vector space of finite singular, piecewise
smooth $k$-chains in a space $X$. Consider the  projection 
$\wb M\to\wh M$ obtained by collapsing each boundary torus to a point.
Given a relative cycle 
$c=\sum_i\lambda_i\sigma_i$ in $C_k(\wb M)$ 
{\em i.e.} a chain $c$ such that $\partial c\in C_{k-1}(\partial \wb M)$, 
we also call $c$ the chain induced on  $C_k(\wh M)$ with $\partial
c\in C_k($ideal points$)$. We call $\wt c$ a lift of $c$ to $\wh{\wt
  M}$, that is $\wt c=\sum_i\lambda_i\wt\sigma_i\in C_k(\wh{\wt M})$
where each $\wt\sigma_i$ is a lift of $\sigma_i$.

\begin{remark}\label{r remA}
  The chain $\wt c$ in general is not a relative cycle.
Nevertheless, since $c$ is a relative cycle, 
assuming $\partial \wt c=\sum_jl_j\eta_j$,
there exists a family $\{\alpha_j\}$ of elements
  of $\pi_1(M)$ such that 
$$\sum_jl_j\cdot\alpha_{j*}(\eta_j)\in C_k(\textrm{ideal points})$$ 
where $\pi_1(M)$ acts on $\wt M$ via deck transformations and
$\alpha_{j*}(\eta_j)$ is the composition of $\alpha_j$ with $\eta_j$. 
\end{remark}

We set $\wb \sigma_i=(D_{\rho})_*(\wt\sigma_i)$ and $\wb c=\sum_i\lambda_i\wb
\sigma_i=(D_{\rho})_*(\wt c)\in C_k(\wb {\M H}^3)$.  We restrict now the
class of simplices we want to use. 
\begin{defi}
We call a $k$-simplex $\sigma:\Delta^k\to\wb{\M H}^3$ {\em admissible}
if for any sub-simplex $\eta$ of
$\sigma$, if the interior of $\eta$ touches $\partial \M H^3$  then
$\eta$ is constant. 
A chain is admissible if its
simplices are admissible.   
\end{defi}

\begin{lemma}\label{l lemB}
  Let $c'=\sum_i\lambda_i\sigma_i'$ be a relative cycle in $C_k(\wb
  M)$. Then there exists a cycle $c=\sum_i\lambda_i\sigma_i$ (with the
  same
  $\lambda_i$'s) such that $[c']=[c]\in H_k(\wb M,\partial \wb M)$ and
  such that $\wb c$ is admissible. 
\end{lemma}
\proof
For any chain $\beta\in C_k(\wb M)$, define span$(\beta)$ as the set
of all the subsimplices of $\beta$ (of any dimension). Given the chain
$c'$, construct $c$ as follows: near $\partial \wb M$ push $c'$ a little
inside $M$, keeping fixed only the simplices of span$(\partial
c')$. This operation can be made via an homotopy, so
$[c]=[c']$. Moreover, the only simplices of $c$ that touch $\partial
\wb M$ are the ones of span$(\partial c)$. Finally, $c$ is admissible
because, if $\wb\sigma_i(x)\in\partial \M H^3$, then from the
definition of pseudo-developing map it follows that $\sigma_i(x)$ is
an ideal point. Thus $x$ lies on a face $F$ of $\sigma_i$ such that the
simplex $\eta=(\sigma_i)|_{F}$ belongs to ${\rm span}(\partial c)$. It follows
that $\wt\eta$ is a constant map and then also $\wb \eta$ is constant.\qed   

We call $\wb C_k(\wb{\M H}^3)$ the vector space of admissible chains. Note
that the boundary operator is well-defined on 
$\oplus_k\wb C_k(\wb{\M H}^3)$ (The boundary of an admissible cycle is
admissible). 

\begin{defi}
  For any admissible simplex $\sigma:\Delta^k\to\wb{\M H}^3$,  a
  straightening ${\rm Str}(\sigma):\Delta^k\to\wb{\M H}^3$ is a simplex  
  agrees with $\sigma$ on the $0$-skeleton, moreover we require
  ${\rm Str}(\sigma)$ to be a standard map whose image is the convex hull of
  its vertices. For any chain $c=\sum_i\lambda_i\sigma_i$ a
  straightening of $c$ is a chain
  ${\rm Str}(c)=\sum_i\lambda_i{\rm Str}(\sigma_i)$.  
\end{defi}

A straightening of a simplex is admissible because any straight
  simplex is admissible. The straightening of a
 simplex is not unique in general.
 Nevertheless, as the following lemma shows, 
 it is possible to choose a straightening
  for any simplex compatibly with the boundary operator of
  $\oplus_k\wb C_k(\wb{\M H}^3)$.

\begin{lemma}
    There exists a chain-map 
${\rm Str} :\oplus_k\wb C_k(\wb{\M H}^3)
\to\oplus_k\wb C_k(\wb{\M H}^3)$ 
that maps
each simplex to one of its straightenings and such that for any
isometry $\gamma$ of $\M H^3$, $\gamma_*\circ{\rm Str}={\rm Str} 
\circ\,\gamma_*$.
 \end{lemma}
\proof
Let $K$ be the set of pairs $\{(B,f)\}$ where 
$B$ is a sub-space of $\oplus_k\wb C_k(\wb{\M H}^3)$ and 
$f:B\to\oplus_k\wb C_k(\wb{\M H}^3)$ is a linear map, such that:
\begin{itemize}
\item $\partial(B)\subset B$.
\item  $\forall\gamma\in{\rm Isom}(\M H^3),\ \gamma_*(B)\subset B$.
\item $\forall \sigma\in B$, $f(\sigma)$ is a straightening of
  $\sigma$.
\item $\forall \gamma\in$ {\rm Isom}$(\M H^3)$,
  $f\circ\gamma_*=\gamma_*\circ f$.
\item $f\circ\partial=\partial\circ f$.
\end{itemize}

Note that $K$ is not empty because each $0$-simplex is admissible and
it is itself its unique straightening, so that $(\wb C_0(\wb{\M
  H}^3),Id)\in K$.
We order $K$ by inclusion ({\em i.e.} $(B,f)\prec(C,g)$ iff $B\subset C$ and
$g|_{B}=f$) and use Zorn's lemma. Let $\{(B_\xi,f_\xi)\}$ be an
ordered sequence in $K$. Clearly 
$$(B_\infty=\cup_{\xi} B_\xi,f_\infty=
\cup_{\xi} f_\xi)$$ 
is an upper bound for $\{(B_\xi,f_\xi)\}$. It follows that there
exists a maximal element $(\wb B,\wb f)\in K$. We claim that 
$\wb B=\oplus_k\wb C_k(\wb{\M H}^3)$. Suppose the contrary. Let 
$k=\min\{n\in\M N: \wb C_n(\wb{\M H}^3)\not\subset\wb B\}$ and let
$\sigma$  be a simplex of $\wb C_k(\wb{\M H}^3)\setminus \wb B$. 
If $k=0$, set $B_1$ the space spanned by 
$\wb B$ and $\displaystyle{\bigcup_{\gamma\in{\rm Isom}(\M H)^3}
  \gamma_*(\sigma)}$, define
$\wb f(\sigma)=\sigma$, $\wb f(\gamma_* (\sigma))=\gamma_*(\wb
f(\sigma))$ and extend $\wb f$ on $B_1$ by linearity.
Then
$$(\wb B,\wb f)\prec(B_1,\wb f)$$
contradicting the maximality of $(\wb B,\wb f)$. 
If $k>0$, then $\wb f$ is defined
on $\partial\sigma$ and, as $\wb f(\partial \sigma)$ is standard, it
is not hard to 
show that it extends to a standard map $\wb f(\sigma)$ 
defined on the whole $\Delta^k$.
Then define $B_1$ and extend $\wb f$ to $B_1$ as above.
 Again we have 
$(\wb B,\wb f)\prec(B_1,\wb f)$, that
contradicts the maximality of $(\wb B,\wb f)$. 

Thus $\wb B=\oplus_k\wb C_k(\wb{\M H}^3)$ and $\wb f$ is the
requested chain map Str.\qed

  \begin{lemma}\label{l lemC}
    There exists a homotopy operator $H:\oplus_k\wb C_k(\wb{\M
  H}^3)\to\oplus_k\wb C_k(\wb{\M H}^3)$ between Str and the
  identity such that  $H\circ\gamma_*=\gamma_*\circ H$ for any
  isometry $\gamma$ of $\M H^3$.  
  \end{lemma}
\proof 
A homotopy operator between ${\rm Str}$ and ${\rm Id}$ 
is a chain map of degree $1$, {\em i.e.} a map
$H:\wb C_k(\wb{\M H}^3)\to \wb C_{k+1}(\wb{\M H}^3)$,  such that 
$$\textrm{Str} -Id =\partial \circ H- H\circ\partial.$$ 
For any admissible $\sigma:\Delta^k\to\wb{\M H}^3$, let $h_\sigma(t,x)$ be
the homotopy constructed as follows: $h_\sigma(t,x)$ is the convex
combination from $\sigma(x)$ to ${\rm Str}(\sigma)(x)$ 
if $\sigma(x)\notin\partial
\M H^3$ and $h_\sigma(t,x)=\sigma(x)$ otherwise. Note that from the
admissibility of $\sigma$ it follows that $h_\sigma(\infty,x)={\rm
  Str}(\sigma)(x)$ for any $x$. 
So the $h_\sigma$ actually is a homotopy between $\sigma$ and ${\rm
  Str}(\sigma)$. 

As $h_\sigma$ is a map $h_\sigma:\Delta^k\times[0,\infty]\to\wb{\M
  H}^3$, up to triangulating $\Delta^k\times[0,\infty]$, it is a
  chain in $C_{k+1}(\wb{\M H}^3)$. Fix a canonical triangulation of
  $\Delta^k\times[0,\infty]$ and define $H(\sigma)$ as $h_\sigma$. Since
  $$\partial(\Delta^k\times[0,\infty])=
\Delta^k\times\{\infty\}- 
\Delta^k\times\{0\}+
\partial\Delta^k\times[0,\infty]$$ then 
$\partial \circ H={\rm Str}- Id+H\circ\partial$.

Since $h_\sigma$ is constructed using geodesic
rays, then for every isometry $\gamma$ we 
have $h_{\gamma\circ\sigma}=\gamma\circ h_\sigma$. It follows that
$H\circ\gamma_*= \gamma_*\circ H$.

Finally, admissibility of $h_\sigma$ follows from admissibility of $\sigma$.
\qed
\begin{lemma}\label{l lemD}
  Let $c=\sum_i\lambda_i\sigma_i$ be a chain in $C_k(\wb
  M)$. Let $\{\gamma_j\}$ be a finite set of isometries
  and let $A$ be the hyperbolic convex hull in $\M H^3$ of 
$$\bigcup_{i,j} \gamma_j({\rm Im}(\wb\sigma_i)).$$
Then $A$ has finite volume.
\end{lemma}
\proof
Since $D_\rho$ has the  cone property on the cusps and since $c$ is a
finite sum of simplices, then $A$ is
contained in a geodesic polyhedron with a finite number of vertices,
and such a polyhedron  has finite volume.\qed 

We are now ready to complete the proof of Theorem~\ref{t sv}.
Let $c=\sum_i\lambda_i\sigma_i$ be a relative cycle in  
$C_3(\wb M)$ such that $[c]=[\wb M]$ in $H_3(\wb M,\partial \wb M)$. 
By Lemma~\ref{l lemB} we can suppose that $c$ is admissible.
Assume
 $\partial \wt c=\sum_jl_j\eta_j$. 
By Remark~\ref{r remA}, there
exists a finite set $\{\alpha_j\}\subset\pi_1(M)$ such that
$\sum_jl_j\cdot\alpha_{j_*}\eta_j\in C_2($ideal points$)$.

Let $A$ be as in Lemma~\ref{l lemD}, where we use 
$\{\rho(\alpha_j)\}\cup\{{\rm Id}\}$ as the set of isometries.
Since $A$ has finite volume, then the volume form $\omega$ of $\M H^3$
is an element of $H^3(A)$. Moreover, the straightening of any
admissible simplex
in $\wb C_k(A)$ is contained in $\wb C_k(A)$ and, since the homotopy operator
$H$ is constructed  using convex combinations, then $H$ is
well-defined on $\oplus_k\wb C_k(A)$. Called $v_i$ the volume of the
straight version of $\sigma_i$, we have
$$
\begin{array}{lcl}
\vol(\rho)&=&\langle D_\rho^*\omega,\wt c\rangle =
\langle\omega,(D_\rho)_*(\wt c)\rangle =
\langle\omega,\wb c\rangle\\
&=& \langle\omega,\textrm{Str}\wb c\rangle+ \langle\omega,H\partial\wb
c\rangle- \langle\omega,\partial H \wb c\rangle\\
&=&
\sum_i\lambda_iv_i+\langle\omega,H\partial \wb c\rangle- \langle d\omega,H\wb
c\rangle=\sum_i\lambda_iv_i+\langle\omega,H\partial \wb c\rangle
\end{array}
$$
By Lemma~\ref{l lemC} we have
$\rho(\alpha_j)_*H=H\rho(\alpha_j)_*$. Moreover, the volume form is
invariant by isometries. 
 It follows that
$$
\begin{array}{lcl}
\langle\omega,H\partial \wb c\rangle&=&
\langle\omega,H\sum_jl_j(D_\rho)_*\eta_j\rangle=
\sum_jl_j\langle\omega,H(D_\rho)_*\eta_j\rangle\\
&=&
\sum_jl_j\langle\rho(\alpha_j)^*\omega,H(D_\rho)_*\eta_j\rangle=
\sum_jl_j\langle \omega,\rho(\alpha_j)_*H(D_\rho)_*\eta_j\rangle\\
&=&
\sum_jl_j\langle \omega,H\rho(\alpha_j)_*(D_\rho)_*\eta_j\rangle
=\sum_jl_j\langle \omega,H(D_\rho)_*\alpha_{j*}\eta_j\rangle\\&=&
\langle \omega,H(D_\rho)_*\sum_jl_j\alpha_{j*}\eta_j\rangle.
\end{array}
$$
 The last product is zero because $D_{\rho*}\sum_jl_j\alpha_{j*}\eta_j$ lies
 on the ideal points of $A$, where $H$ is fixed and
 $\omega$ vanishes.

This completes the proof of Theorem~\ref{t sv}.\qed

\begin{remark}\label{r difp}
  Theorem~\ref{t sv} applies when the cycle $c$ is an ideal
  triangulation. So it implies Theorem~\ref{t vsd}. 
\end{remark}

\begin{cor}\label{c vg}
  Let $M$ be a graph $3$-manifold. Then for all representations
  $\rho:\pi_1(M)\to{\rm Isom}^+(\M H^3)$ we have $\vol(\rho)=0$.
\end{cor}
\proof This is because for each graph manifold $M$ we have
$||(M,\partial M)||=0$ (see~\cite{G}, \cite{Ku}).\qed
\begin{cor}\label{c df}
  Let $M$ be a complete hyperbolic $3$-manifold of finite volume. Then for all
  representations $\rho:\pi_1(M)\to{\rm Isom}^+(\M H^3)$ we have
  $|\vol(\rho)|\leq$ $\vol(M)$.
\end{cor}
\proof This follows from the fact that for complete hyperbolic
$3$-manifold we have $\vol(M)=v_3||(M,\partial M)||$ (see~\cite{G},
\cite{Ku}).\qed 

In \cite{D} it is proved that, for compact manifolds, 
equality holds if and only if
$\rho$ is discrete and faithful. 
In the next section we show that this is true in general for manifolds
of finite volume.

\section{Rigidity of representations}\label{s r}
This section is completely devoted to proving the following:
\begin{teo}\label{t ug}
  Let M be a non-compact, complete, orientable hyperbolic $3$-manifold
  of finite 
  volume. Let $\Gamma\cong\pi_1(M)$ be the sub-group of ${\rm PSL}(2,\M C)$
  such that $M=\M H^3/\Gamma$. Let 
  $\rho:\Gamma\to$ {\rm PSL}$(2,\M C)$ be a representation.
If $|\vol(\rho)|=\vol(M)$ then $\rho$ is discrete and
faithful. More precisely there exists
$\f\in{\rm PSL}(2,\M C)$ such
that for any $\gamma\in\Gamma$ 
$$\rho(\gamma)=\f\circ\gamma\circ\f^{-1}.$$
\end{teo}

\begin{remark} It is
  well-known that, in the hypotheses of Theorem~\ref{t ug}, the
  manifold $M$ is the interior of a compact
  manifold $\wb M$ whose boundary consists of tori. Thus $M$ is a
  cusped manifold and, 
  by Remark~\ref{r acmcbid}, all the definitions and results we gave
  for ideally triangulated manifold apply. 

As product structure on the cusps
  we fix the horospherical one, having the arc-length as cone
  parameter. 
For this section $D_\rho$ will denote a fixed pseudo-developing map 
for $\rho$.
\end{remark}

\begin{remark}
By Proposition~\ref{p pos} we can suppose $\vol(\rho)\geq0$. 
\end{remark}

\begin{remark}\label{r none}
  A subgroup of {\rm PSL}$(2,\M C)$ is said to be {\em elementary} if it has an
  invariant set of at most two points in $\partial \M H^3$. If the
  image of $\rho$ is elementary, then one can construct a
  pseudo-developing map as in Lemma~\ref{l 1} in such a way that 
  all the tetrahedra of any ideal
  triangulation of $M$ collapse in the straightening. 
Thus, by Theorem~\ref{t vsd}, $\vol(\rho)=0$. 
\end{remark}
This remark implies that, in the present case, 
since $\vol(\rho)=\vol(M)\neq0$, the image of
$\rho$ is non-elementary.

The idea for proving Theorem~\ref{t ug} is to
rewrite the Gromov-Thurston-Goldman-Dunfield proof of Mostow's
rigidity, valid in the compact case.

We will follow the lead-line of~\cite{D}, 
with the difference that we will use classical chains
instead of measure-chains. The technique for constructing classical chains
representing smear-cycles is that used in~\cite{BP} for the proof
of Mostow's rigidity for compact manifolds. As an effect of
non-compactness we will work with infinite chains. Therefore, we have to
prove that some usual homological arguments actually work for our
chains.

The core of the proof is to deduce from the equality
 $\vol(\rho)=\vol(M)$ that $D_\rho$
``does not shrink the volume.'' This allows us to construct 
a measurable extension of $D_\rho$ to the whole $\wb{\M H}^3$,
whose restriction to $\partial \M H^3$ is almost everywhere
a M\"obius transformation. 
Such a M\"obius transformation will be the $\f$ of Theorem~\ref{t ug}. 

The key fact is the following proposition,  whose proof can be found 
in~\cite{D}~(claim $3$ of Theorem 6.1).
\begin{prop}\label{p k}
Let $f:\partial \M H^3\to \partial \M H^3$ be a measurable 
map that maps the vertices of almost all regular ideal tetrahedra to
vertices of regular ideal tetrahedra. 
Then $f$ coincides almost everywhere with the trace of an isometry \f. 
\end{prop}

We want to apply Proposition~\ref{p k} to $D_\rho$, 
and we do it in two steps. 
Let $M_0$ be $M$ minus the cusps and let $\pi:\M H^3\to M$ be the
universal cover.

\begin{prop}\label{p nnp1}
The map $D_\rho$ extends to $\wb{\M H}^3$. More precisely, there exists 
a measurable map $\wb D_\rho:\partial \M
H^3\to\partial \M H^3$ such that for
almost all $x\in\partial \M H^3$, for any geodesic $\gamma^x$
ending at $x$, for any sequence $t_n\to\infty$ such that
$\pi(\gamma^x(t_n))\in M_0$, we have
$$\lim_{n\to\infty}D_\rho(\gamma^x(t_n))=\wb D_\rho(x).$$ 
\end{prop}

\begin{prop}\label{p nnp2}
The map $\wb D_\rho$ satisfies the hypothesis of Proposition~\ref{p k}. 
\end{prop}

Before proving Propositions~\ref{p nnp1}, and~\ref{p nnp2} we show how
they imply Theorem~\ref{t ug}  

\noindent{\em Proof of \ref{t ug}.}
By Proposition~\ref{p nnp2}, Proposition~\ref{p k} applies.
By Proposition~\ref{p nnp1}  
the equivariance of $D_\rho$ implies the equivariance of $\wb D_\rho$, 
getting for any $\gamma\in\Gamma$
$$\rho(\gamma)=\f\circ\gamma\circ\f^{-1}.$$\qed
\begin{remark}
Both Propositions~\ref{p nnp1}, and~\ref{p nnp2} will follow from  
Lemma~\ref{l cind} and~\ref{l ladf} below.
We notice that Lemma~\ref{l cind} is a restatement of Lemmas~$6.2$
of~\cite{D}. While Proposition~\ref{p nnp2} corresponds to {\em Claim
  $2$} of~\cite{D}.
Proposition~\ref{p nnp1} follows from 
Lemmas~\ref{l cind} and~\ref{l ladf} exactly as in~\cite{D}. 
We will give a complete proof of Proposition~\ref{p nnp2} because 
the proof of {\em Claim $2$} in~\cite{D} seems to be incomplete.
\end{remark}

From now until Lemma~\ref{l nnp3} we will  describe how 
to construct a simplicial version of the 
{\em smearing} process
of  measure-homology (see~\cite{Th1} or~\cite{R}).
Then we will prove Lemma~\ref{l ladf}. Finally we will complete the proof
of Theorem~\ref{t ug} by proving Propositions~\ref{p nnp1} 
and~\ref{p nnp2}.

Let $\mu$ be the Haar measure on Isom$(\M H^3)$ such that for each
$x\in\M H^3$ and $A\subset \M H^3$ we have
$$\mu\{g\in{\rm Isom}(\M H^3)\, :\, g(x)\in A\}=\vol(A)$$
where $\vol(A)$ is the hyperbolic volume of $A$.

In the following by a tetrahedron of $\wb{\M H}^3$
we mean an ordered 4-tuple of points (the vertices). 
The volume of a tetrahedron is the hyperbolic volume with sign of
the convex hull of its vertices. 

Let $S$ be the set of all genuine (non-ideal, non-degenerate) tetrahedra:
$$S=\{(y_0,\dots,y_3)\in(\M H^3)^4:\ \vol(y_0,\dots,y_3)\neq0\}.$$

For any $Y\in S$ let $S(R)$ be the set of all isometric copies of $Y$:
$$S(Y)=\{X\in S:\ \exists g\in{\rm Isom}(\M H^3),\  X=g(Y)\}.$$

Then a natural bijection 
$f_Y:\ ${\rm Isom}$(\M H^3)\to S(Y)$ is well-defined by
$$f_Y(g)=g(Y).$$
Thus $\mu$
induces a measure, which we still call $\mu$, on $S(Y)$ defined by
$$\mu(A)=\mu(f_Y^{-1}(A)).$$

We consider the sets  $S_\pm(Y)=f_Y^{-1}(${\rm Isom}$^\pm(\M H^3))$  
of tetrahedra respectively positively and negatively isometric to $Y$. Note
 that $S_+(Y)$ and $S(Y)_-$ are both measurable.

Set $\mfk S=\Gamma^4/\Gamma$ where $\Gamma$ acts on $\Gamma^4$ by
left multiplication. 
Each element $\sigma=[(\gamma_0,\dots,\gamma_3)]\in \mfk S$ 
has a unique representative with $\gamma_0={\rm Id}$. When we
write $\sigma\in\mfk S$ we tacitly assume that the representative of
the form $(\gamma_0,\dots,\gamma_3)$ with $\gamma_0={\rm Id}$ has been chosen.
So $\gamma_0$ is always the identity.

For the rest of the section we fix a fundamental polyhedron 
$F\subset\M H^3$ for $M$.
For all $\e>0$ let $\call F^\e$ be a locally finite \e-net in $F$.
For any $\xi\in\call F^\e$ let 
$$F_\xi=\{x\in F:\ d(x,\xi)=d(x,{\call F}^\e)\}.$$

Each $F_\xi$ is a geodesic polyhedron of diameter less than
\e. From the cone-property of $D_\rho$ it follows
that the diameters of $D_\rho(F_\xi)$ are bounded by a constant
$\delta$
that depends on \e.
Moreover, by removing some boundary face from some $F_\xi$, we get
that $F$ is the disjoint union of the $F_\xi$'s. We define now a
family of special simplices. Let
$$\mfk N=\{(\gamma_0,\dots,\gamma_3,\xi_0,\dots,\xi_3):\
(\gamma_0,\dots,\gamma_3)\in{\mfk S},\  
\xi_i\in{\call F}^\e \textrm{ for all } i\}.$$ 

For each $\eta\in\mfk N$ define $\Delta_\eta$ as the 
straight geodesic singular $3$-simplex whose vertices are the  points
$\xi_0,\gamma_1(\xi_1),\gamma_2(\xi_2),\gamma_3(\xi_3)$, more precisely
$$\Delta_\eta:\Delta^3\ni t \mapsto 
\pi\Big(\sum_{i=0}^3t_i\gamma_i(\xi_i)\Big).$$

For each tetrahedron $X=(x_0,\dots,x_3)\in S(Y)$ 
there exists a unique 
$\eta=(\gamma_0,\dots,\gamma_3,\xi_0,\dots,\xi_3)\in\mfk N$ such that
$x_i\in \gamma_i(F_{\xi_i})$ for $i=0,\dots,3$. This defines a function
$$s_Y:S(Y)\to\mfk N.$$
Roughly speaking, \mfk N is a locally finite \e-net in the space of
$3$-simplices of $M$ and $s_Y$ is the ``closest point''-projection.

For any $\eta\in\mfk N$ define
$$a_Y^\pm(\eta)=\mu\{s_Y^{-1}(\eta)\cap S_\pm(Y)\}=
\mu\{X\in S_\pm(Y) :\ x_i\in\gamma_i(F_{\xi_i})\}$$
and $$a_Y(\eta)=a^+_Y(\eta)-a^-_Y(\eta).$$

In the language of measures, one can think of $a_Y^\pm$ 
as the push-forward of the
measure $\mu$ under the map $s_Y:S_\pm(Y)\to\mfk N$.
This is the key for the passage from measure-chains to classical ones.

The smearing of the tetrahedron $Y$ is the cycle:
$$Z_Y=\sum_{\eta\in\mfk N}a_Y(\eta)\Delta_\eta.$$
We notice that, as $\mfk N$ depends on the family $\call F^\e$, 
the cycle $Z_Y$ actually depends on \e.

\begin{remark}
  The smearing of a tetrahedron in general is not a finite
  sum. Nevertheless, as the following lemma shows, it has bounded $l^1$-norm.
\end{remark}

\begin{lemma}\label{l L1} For any $Y\in S$, we have
  $\displaystyle{\sum_\eta|a_Y(\eta)|<\vol(M)}.$
\end{lemma}
\proof If $Y=(y_0,\dots y_3)$ then
$$
\begin{array}{l}
\displaystyle{\sum_\eta|a_Y(\eta)|\leq\sum_\eta
  \big(a^+_Y(\eta)+a_Y^-(\eta)\big)}=\displaystyle{
  \sum_\eta\mu\{s_Y^{-1}(\eta)\}=\mu\Big\{\,\bigcup_\eta
  s_Y^{-1}(\eta)\Big\}}\\
=\displaystyle{\mu\{\,s_Y^{-1}(\mfk
  N)\,\}}=
\mu\{\,f_Y^{-1}s_Y^{-1}(\mfk N)\,\}=
 \mu\{g:\ g(y_0)\in F\}=\vol(F)=\vol(M).
\end{array}
$$
\qed
\begin{lemma}\label{l nnp3}
  The infinite chain $Z_Y$ is a cycle, {\em i.e.} $\partial Z_Y=0$.
\end{lemma}
\proof
First note that the $l^1$-norm of $\partial Z_Y$ is bounded by
$4$ times the $l^1$-norm of $Z_Y$. Thus all the sums we will consider
make sense.

Let $\upsilon$ be a simplex of $\partial Z_Y$. By construction
$\upsilon$ is obtained as the projection of an $(n-1)$-simplex having
vertices in $F_{\xi_0},\gamma_1(F_{\xi_1}),\gamma_2(F_{\xi_2})$ for
some $\gamma_1,\gamma_2\in\Gamma$ and $\xi_0,\xi_1,\xi_2\in\call
F^\e$. Let $A_\upsilon$ be the set of the elements of \mfk N of the
form
$\eta=(\gamma_0,\gamma_1,\gamma_2,\gamma,\xi_0,\xi_1,\xi_2,\xi)$ with
$\gamma\in\Gamma$ and $\xi\in\call F^\e$.
The simplices $\Delta_\eta$ of $Z_Y$ having $\upsilon$ as the last face
contribute to the coefficient of $\upsilon$ in $\partial Z_Y$ by
$$\begin{array}{ll}
\displaystyle{\sum_{\eta\in A_\upsilon} a_Y(\eta)=\sum_{\eta\in
  A_\upsilon}\mu(s_Y^{-1}(\eta)\cap S_+(Y)) -\sum_{\eta\in A_\upsilon}
  \mu(s_Y^{-1}(\eta)\cap S_-(Y))}\\=
\mu(s_Y^{-1}(A_\upsilon)\cap S_+(Y))-\mu(s_Y^{-1}(A_\upsilon)\cap
S_-(Y))=0.\end{array}$$ 

The same calculation, made with the simplices having $\upsilon$ as the
$i$th face, shows that the coefficient of $\upsilon$ in $\partial Z_Y$
is zero.\qed

For any ideal, non-flat, tetrahedron $Y=(y_0,\dots,y_3)$ 
let $t\mapsto y_i(t)$ be the geodesic ray from the center of mass of $Y$ 
to $y_i$, $i=0,\dots,3$.
For any $R>0$ let $Y_R$ be the following element of $S$:
 $$Y_R=(y_0(R),\dots,y_3(R)).$$ 

\begin{remark}
From now on we fix a positively oriented regular ideal tetrahedron $Y$, 
and we write $S_\pm(R)$, $f_R$, $s_R$,
$a_R(\eta)$ and $Z_R$
for $S_\pm(Y_R)$, $f_{Y_R}$, $s_{Y_R}$, $a_{Y_R}(\eta)$ and $Z_{Y_R}$.
\end{remark}

We say that a $3$-simplex $\Delta$ is \e-close to a tetrahedron $X$ if the
vertices of $\Delta$ are \e-close to $X$. We define
$$\epsilon(R,\e)=\sup\{v_3-\vol(\Delta):\ 
\Delta\textrm{ is \e-close to an element of } S(R)\}$$

  \begin{lemma}\label{l and}
    For any fixed \e,  for large $R$ the function $\epsilon(R,\e)$
    goes to zero exponentially in $R$.
  \end{lemma}
This is because $v_3-\vol(Y_R)$ goes to zero like $e^{-R}$ and
the volume of any $\Delta$ which is \e-close to $Y_R$ is close to the volume
of $Y_R$. See~\cite{BP},~\cite{D},~\cite{Th1} for details.

\begin{remark}
What we actually need to prove our claims is a restatement for $Z_R$ 
of the {\em Step~\ref{stp3}} of Theorem~\ref{t sv}. 
From now until Proposition~\ref{p 27}, we prove facts that are
standard for finite chains, but need a proof for $Z_R$.
\end{remark}

For $\eta\in\mfk N$, we set $v_\eta=\vol(\Delta_\eta)$.
Using the fact that all the $F_\xi$'s have diameter less than \e, one
can prove the following lemma (see~\cite{BP} for details). Recall that
$\mfk N$ depends on $\call F^\e$ and so it depends on \e.
\begin{lemma}\label{l novo}
  For any $\e>0$, for large enough $R$  
  we have that for any $\eta\in \mfk N$ 
  \begin{itemize}
  \item $a_R^+(\eta)\cdot a_R^-(\eta)=0$.
  \item $a_R(\eta)\neq 0\implies a_R(\eta)\cdot v_\eta\geq 0$. 
  \end{itemize}
\end{lemma}

\begin{lemma}\label{l lemmaX}
There exists a constant $c$ such that $|D_\rho^*\omega|<c|\omega|$,
where $\omega$ is the volume-form of $\M H^3$.
\end{lemma}
\proof
Let $M_0$ be $M$ minus the cusps.
The function $|D_\rho^*\omega|/|\omega|$ is
continuous and hence bounded on $M_0$. In the cusps, by direct
calculation and  using the cone property of $D_\rho$, one can show
that the same bound holds.\qed

\begin{lemma}
The integrals $\langle\omega,Z_R\rangle$ and $\langle
D_\rho^*\omega,Z_R\rangle$ are well-defined.
\end{lemma}
\proof
As $\sum|a_R(\eta)|<+\infty$, since
$|\langle\omega,\Delta_\eta\rangle|$ is bounded by $v_3$, then $\langle
\omega,Z_R\rangle$ is well-defined. Consider now $D_\rho^*\omega$. 
From Lemma~\ref{l lemmaX} it follows that the integral of $|D_\rho^*|$
over straight geodesic simplices is bounded by $cv_3$. Hence also
$\langle D_\rho^*\omega,Z_R\rangle$  is well-defined.\qed

As above, let $M_0$ denote $M$ minus the cusps and, 
for $k\in \M N^*$ let 
$$M_k=\displaystyle{\bigcup_{T\subset\partial M_0}T\times[k-1,k).}$$  

Let ${\call F}_k^\e={\call F^\e}\cap \pi^{-1}(M_k)$ and 
$\mfk N_k=\{\eta\in\mfk N:\ \xi_0\in{\cal F}_k^\e\}$.
We have 
$$Z_R=\sum_{k\in\M N}\sum_{\eta\in\mfk N_k}a_R(\eta)\Delta_\eta.$$

\begin{lemma}\label{l lemmaY}
  For any $k$ the chain $\sum_{\eta\in\mfk N_k} a_R(\eta)\Delta_\eta$ is a
  finite sum.
\end{lemma}
\proof
If $a_R(\eta)\neq0$ and $\eta\in\mfk N_k$ then $\Delta_\eta$ is
\e-close
to an element $X\in
S(R)$ having
first vertex in $F_{\xi_0}$ with $\xi_0\in{\cal F}^\e_k$ . 
Since $\call F^\e$ is
locally finite and $\wb M_k$ is compact,  ${\cal F}_k^\e$ is
finite, so there is only a finite number of possibilities for $\xi_0$.
Since $\wb F_{\xi_0}$ is compact, 
 any $X\in S(R)$
with first vertex in $F_{\xi_0}$ lies on a compact ball $B$ of $\M H^3$.
Since $F$ is a fundamental domain, then  
there exists only a finite number of elements $\gamma\in\Gamma$ so
that $\gamma(F)$ intersects $B$. Then for any $\xi_0$ there is only a
finite number of possibilities for $\xi_1,\ \xi_2$ and $\xi_3$. 
It follows that there exists only a finite number of
$\eta\in\mfk N_k$ such that $a_R(\eta)\neq0$.\qed

\begin{lemma}\label{l lemmaZ}
For any $R$, if $k$ is large enough, then for any $\eta\in \mfk N_k$
  with $a_R(\eta)\neq 0$, the simplex $\Delta_\eta$  
  is completely contained in a cusp of $M$. 
\end{lemma}
\proof
If $X=(x_0,\dots,x_3)\in S(R)$ then $X$ lies in the ball $B(x_0,2R)$. 
Since $M$ has a finite number of cusps, for any $R$ there
exists $m\in \M N$ such that for $k\geq m$ if $x_0\in M_k$ 
then the whole ball
$B(x_0,2R+\e)$ is contained in the cusp containing $x_0$.
If $\eta\in\mfk N_k$ and $a_R(\eta)\neq0$, then there exists
$X\in S(R)$ with $x_0\in \pi^{-1}(M_k)\cap F$ 
hence $\Delta_\eta$ is \e-close to $X$. Thus $\Delta_\eta\subset
B(x_0,2R+\e)$ is contained in the cusp that contains $x_0$.  
\qed 

Now for $k\in\M N$ define
$$Z_{R,k}=\sum_{j<k}\sum_{\eta\in \mfk N_j}a_R(\eta)\Delta_\eta.$$

$Z_{R,k}$ is a finite chain by Lemma~\ref{l lemmaY}. 
Moreover, since $\partial Z_R=0$, 
then each simplex $\upsilon$ 
of $\partial Z_{R,k}$ appears as a face of a simplex $\Delta_\eta$
with $a_R(\eta)\neq 0$ and $\eta\in{\mfk N}_j$ for some  $j\geq k$. Therefore,
by Lemma~\ref{l lemmaZ}, for $k$ large enough each simplex $\upsilon$ 
of $\partial Z_{R,k}$ is contained in a cusp of $M$. Thus to each $\upsilon$
 there corresponds an ideal point of $\wh M$. For each $\upsilon\in\partial
  Z_{R,k}$ let $\lambda_{R,k}(\upsilon)$ be the coefficient of $\upsilon$ in
  $\partial Z_{R,k}$ and let $C_\upsilon$ be the cone from $\upsilon$
to the corresponding ideal point.

Let $\wh Z_{R,k}$ be the chain obtained by adding to $Z_{R,k}$ the
cones $C_\upsilon$:
$$\wh Z_{R,k}=Z_{R,k}+\sum_{\upsilon\in\partial Z_{R,k}}\lambda_{R,k}(\upsilon)C_\upsilon.$$
  The chain $\wh Z_{R,k}$ is a finite sum and it is easily
checked that it is a cycle.

For any $3$-simplex $\Delta$ let $\svol(\Delta)$ denote the volume of
the convex hull of the vertices of $D_\rho(\Delta)$.
For any $\eta\in\mfk N$ set $w_\eta=\svol(\Delta_\eta)$. 

\begin{prop}\label{p 27} For any $R>0$ We have:
  $$\sum_\eta a_R(\eta) v_\eta=
\langle \omega,Z_R\rangle=\langle
D_\rho^*\omega,Z_R\rangle=\sum_\eta a_R(\eta)w_\eta.$$
\end{prop}
\proof The first equality is tautological. We use now the cycles $\wh
Z_{R,k}$ to approximate $Z_R$.
Since $\vol(\rho)= \vol(M)$, then $[\omega]=[D_\rho^*\omega]$ as
elements of $H^3(\wh M)$. Thus for any $k\in\M N$ we have
$\langle \omega,\wh Z_{R,k}\rangle=\langle
D^*_\rho\omega,\wh Z_{R,k}\rangle.$
As in {\em Step~\ref{stp3}} of Theorem~\ref{t sv}, we can straighten the
finite cycle $\wh Z_{R,k}$, getting:
$$\langle \omega,\wh Z_{R,k}\rangle=\langle
D^*_\rho\omega,\wh Z_{R,k}\rangle=
\sum_{j<k}\sum_{\eta\in\mfk N_k}a_R(\eta)w_\eta
+\sum_{\upsilon \in\partial Z_{R,k}}\lambda_{R,k}(\upsilon)\svol(C_\upsilon).$$

For each simplex $\alpha$ of $\wh Z_{R,k}$ we have
$|\vol(\alpha)|\leq v_3$, 
$|\svol(\alpha)|\leq v_3$ and, 
by Lemma~\ref{l lemmaX}, $|\langle D^*_\rho \omega,
\alpha\rangle|\leq cv_3$.  It follows that to get the remaining
inequalities it suffices to show that
$$\lim_{k\to \infty}\sum_{\upsilon\in\partial 
Z_{R,k}}|\lambda_{R,k}(\upsilon)|=0.$$

Since $\partial Z_R=0$, if $\upsilon\in\partial Z_{R,k}$ then
$\upsilon\in\partial \Delta_\eta$ with 
$a_R(\eta)\neq 0$  and $\eta\in\mfk N_j$ for some $j\geq k$. 
So we have
$$
\begin{array}{cl}
~&\displaystyle{\sum_{\upsilon\in\partial Z_{R,k}}|\lambda_{R,k}(\upsilon)|\leq
4\sum_{j\geq k}\sum_{\eta\in\mfk N_j}a^+_R(\eta)+a^-_R(\eta)=
4\sum_{j\geq k}\sum_{\eta\in\mfk N_j}\mu\{s_R^{-1}(\eta)\}}\\
=& \displaystyle{4\sum_{j\geq k} \mu\{s_R^{-1}(\mfk
N_j)\}=4\sum_{j\geq k}\mu\{Y\in S(R):\exists\xi\in{\call F}_j^\e,\
y_0\in F_\xi\}}\\
=& \displaystyle{4\sum_{j\geq k}\sum_{\xi\in{\call F}^\e_j}\vol(F_\xi)
\leq4\vol(\bigcup_{j\geq k-\e}M_j)}
\end{array}
$$
The last term goes to zero as $k\to\infty$ because $M$ has finite
volume and the desired equality follows.\qed

Now that we have Proposition~\ref{p 27}, forget about the cycles 
$\wh Z_{R,k}$.

From triangular inequality, Proposition~\ref{p 27} and 
Lemma~\ref{l novo} we have 
$$
\begin{array}{cl}
&\displaystyle{\sum_\eta|a_R(\eta)|\cdot|w_\eta|\geq 
\Big|\sum_\eta a_R(\eta)w_\eta\Big|=\Big|\sum_\eta a_R(\eta)v_\eta\Big|}\\
=&\displaystyle{
\sum_\eta|a_R(\eta)|\cdot|v_\eta|\geq \sum_\eta
|a_R(\eta)|(v_3-\epsilon(R,\e))}
\end{array}
$$
from which and Lemma~\ref{l L1} we get:
\begin{prop}\label{p p14}
For $R$ large enough we have
$$\sum_{\eta\in\mfk N}|a_R(\eta)|(v_3-|w_\eta|)\leq
\sum_{\eta\in\mfk N}|a_R(\eta)|\epsilon(R,\e)\leq \vol(M)\epsilon(R,\e).$$
\end{prop}

For any $R>0$ let $A_R\subset\mfk N$ be the set of tetrahedra with
``small'' straight volume:
$$A_R=\{\eta\in\mfk N:\ v_3-|w_\eta|>R^2\cdot\vol(M)\cdot\epsilon(R,\e)\}.$$

\begin{lemma}\label{l cind}
  For $R$ large enough we have
$$\sum_{\eta\in A_R}|a_R(\eta)|\leq\frac{1}{R^2}.$$
\end{lemma}
\proof
From Proposition~\ref{p p14} we get
$$
\begin{array}{ll}
&\displaystyle{
R^2\vol(M)\epsilon(R,\e)\cdot\sum_{\eta\in A_R}|a_R(\eta)|
\leq\sum_{\eta\in A_R}|a_R(\eta)|(v_3-|w_\eta|)}\\
\leq&
\displaystyle{\sum_{\eta\in \mfk N}|a_R(\eta)|(v_3-|w_\eta|)}
\leq \vol(M)\epsilon(R,\e)
\end{array}$$
The claimed inequality follows.\qed

\begin{lemma}\label{l ladf}
For almost all isometries $g$ we have
$$\lim_{n\to\infty}\svol(g(Y_n))=v_3.$$  
\end{lemma}

\proof
Since $a_R^+\cdot a_R^-=0$, then 
$\sum_{\eta\in A_R}|a_R(\eta)|=\mu(s_R^{-1}(A_R))$. Thus
for any fixed $R>0$ we have
$$\mu\Big(\bigcup_{\M N\ni n> R}s_R^{-1}(A_n)\Big)
\leq\sum_{n>R}\frac{1}{n^2}<\frac{1}{R}.$$ 

Recalling that for any set $A\subset\mfk N$ we have 
$\mu(s_R^{-1}(A))=\mu(f_R^{-1}s_R^{-1}(A))$, we get
$$\mu\{
g\in\textrm{ Isom}(\M H^3): \exists n>R,\
w_{s_n(g(Y_n))}<v_3-n^2\cdot\vol(M)\cdot\epsilon(n,\e)
\}<\frac{1}{R}.$$ 

From Lemma~\ref{l and} it follows that 
$\lim_{n\to\infty}n^2\epsilon(n,\e)=0$. 
As $R\to\infty$, this implies that for any $\e>0$, for almost any
isometry $g$ we have
$$\lim_{n\to\infty}w_{s_n(g(Y_n))}=v_3.$$

Let $g$ be one of such maps. 
Since the diameters of the $D_\rho(F_\xi)$ are bounded by $\delta$,
then 
$D_\rho(\Delta_{s_R(g(Y_R))})$ is $\delta$-close to
$D_\rho(g(Y_R))$. Recalling that
$w_{s_R(g(Y_R))}=\svol(\Delta_{s_R(g(Y_R))})$, we have that  
$$\lim_{n\to\infty}\svol(\Delta_{s_n(g(Y_n))})=v_3$$
and, since $D_\rho(g(Y_R))$ is $\delta$-close to 
$D_\rho(\Delta_{s_R(g(Y_R))})$, then also
$$\lim_{n\to\infty}\svol(g(Y_n))=v_3.$$
\qed

We sketch here the proof of Proposition~\ref{p nnp1}, referring to~\cite{D}
for details.

\noindent{\em Proof of Proposition~\ref{p nnp1}.}
In the disc model let $\gamma$ be a geodesic from $0$ to a point in
$\partial \M H^3$. Let $X_R$ be a family of regular tetrahedra of edge
$R$ with first vertex in $0$ and second in $\gamma(R)$.
All the claims from Lemma~\ref{l L1} to Lemma~\ref{l ladf} hold for
$\{X_R\}$. It follows that for almost all isometries $g$ we have
$$\lim_{n\to\infty}\svol(g(X_n))=v_3.$$
Then $D_\rho(g(\gamma(n)))$ must reach the boundary of $\M H^3$. 
Using again the above property of the limit,
one can estimate the angle $\alpha(n)$ between the
geodesic from $D_\rho(g(0))$ to $D_\rho(g(\gamma(n)))$ and the geodesic from
$D_\rho(g(0))$ and $D_\rho(g(\gamma(n+1)))$. Such estimate shows that $\sum
\alpha(n)<\infty$, which implies that $D_\rho(g(\gamma(n)))$ converges. 
The claim follows because $D_\rho$ is locally Lipschitz outside the
cusps. Measurability follows because the extension can be viewed as a
point-wise limit of measurable functions.\qed

\begin{remark}
  In general $D_\rho$ is not uniformly continuous in the cusps. So
  it cannot be locally Lipschitz on the whole $\M H^3$.
\end{remark}

We come now to the proof of Proposition~\ref{p nnp2}. 

\begin{lemma}\label{l nnpp}
  Let $X=(x_0,x_1,x_2,x_3)$ be an ideal tetrahedron in $\wb{\M H}^3$. 
Suppose that no three vertices  of $X$ coincide. Then for any $\e>0$ there
exist neighborhoods $U_i$ of $x_i$ in $\wb{\M H}^3$ such that for any
tetrahedron $Y=(y_0,\dots,y_3)$ with $y_i\in U_i$ we have
$|\vol(Y)-\vol(X)|<\e$. 
\end{lemma}

This follows from the formula of the volume for ideal tetrahedra,
see~\cite{BP} for details.

\begin{remark}
  Lemma~\ref{l nnpp} does not hold if three vertices of $X$
  coincide. To see this, let $Y$ be a regular ideal tetrahedron and let
  $\gamma$ be a parabolic or hyperbolic isometry. Then $\gamma^n(Y)$
  is a family of tetrahedra with maximal volume, but at least three of the
  vertices of $\gamma^n(Y)$ converge to the same point.
\end{remark}

\begin{lemma}\label{l era1}
  For almost all regular ideal tetrahedra $Y$, the ideal tetrahedron 
 $\wb D_\rho(Y)$ is defined. Moreover, for almost all $Y$ either 
$\wb D_\rho(Y)$ is
  regular (whence $\vol(\wb D_\rho(Y))=v_3$) or at least three of its vertices
  coincide (whence $\vol(\wb D_\rho(Y))=0$).
\end{lemma}
\proof
Without loss of generality,  we can restrict the first claim 
to the space of positive regular ideal tetrahedra.
We parametrize such a space  with 
$$\{(a,b,c)\in S_\infty^2\times S_\infty^2\times S_\infty^2:\ a\neq b\neq
  c\}$$
where $S_\infty^2=\partial \M H^3$, by mapping $(a,b,c)$ to the unique
positive regular ideal tetrahedron with $(a,b,c)$ as the first three
  vertices. We denote by $Q(a,b,c)$ the fourth vertex of such
  tetrahedron. Since $\wb D_\rho$ is defined almost everywhere,
the first claim follows from Fubini's theorem. The second claim
follows from Lemmas~\ref{l ladf} and~\ref{l nnpp}.\qed

With the above notation, by Lemma~\ref{l era1} we can restate 
Proposition~\ref{p nnp2} as follows. 
\begin{prop}\label{p base}
    The set
$\{Y\in S_\infty^2\times S_\infty^2\times S_\infty^2:\ \vol(\wb
D_\rho(Y))=0\}$ has zero measure.
  \end{prop}
The proof of this result will follow from the next:

\begin{lemma}\label{l basedue}
    If the set
$$\{Y\in S_\infty^2\times S_\infty^2\times S_\infty^2:\ \vol(\wb
D_\rho(Y))=0\}$$ has positive measure, then the map $\wb D_\rho$ is constant
almost everywhere.
\end{lemma}

Before proving Lemma~\ref{l basedue} we show how it
implies Proposition~\ref{p base}.

\ 

\noindent{\em Proof of \ref{p base}.} By contradiction, we apply 
Lemma~\ref{l basedue} deducing that $\wb D_\rho$ is a.e. a constant
$p$.
From the equivariance of $\wb D_\rho$ 
it follows that for any $\gamma \in\Gamma$
and $x\in\partial \M H^3$ we have
$$p=\wb D_\rho \gamma(x)=\rho(\gamma)(\wb D_\rho(x))=\rho(\gamma)(p).$$ 
Thus $p$ is a fixed point of any element of $\Gamma$. This implies that
the image of $\rho$ is elementary, but this cannot happen 
because of  Remark~\ref{r none}.
\qed
  
We now prove Lemma~\ref{l basedue}.
\begin{lemma} In the hypothesis of Lemma~\ref{l basedue}
  there exists a positive-mea\-sure set $A\subset S_\infty^2$ such that
  $\wb D_\rho$ is constant on A. 
\end{lemma}

\proof By Lemma~\ref{l era1} it is not restrictive to suppose that the
set 
$$\{(a,b,c)\in S_\infty^2\times S_\infty^2\times S_\infty^2:\ 
\wb D_\rho(a)=\wb D_\rho(b)=\wb D_\rho(c)\}$$ has positive measure. 
Then by Fubini's theorem there exists a positive-measure set
$A_0\subset S_\infty^2$ such that for all $a_0\in A_0$ the set
$$\{(b,c)\in S_\infty^2\times S_\infty^2:\ 
\wb D_\rho(a_0)=\wb D_\rho(b)=\wb D_\rho(c)\}$$ 
has positive measure in 
$S_\infty^2\times S_\infty^2$. Again by Fubini's theorem for all
$a_0\in A_0$ there exists a positive-measure set $A_1\in S_\infty^2$
such that for any $a_1\in A_1$ the set
$$\{c\in S_\infty^2:\wb D_\rho(a_0)=\wb D_\rho(a_1)=\wb D_\rho(c)\}$$
has positive measure. In particular $\wb D_\rho$ is constant on $A_1$.\qed

We set $p=\wb D_\rho(A_1)$ and $A=\wb D_\rho^{-1}(p)$. 
\begin{remark}
  In the sequel we use the symbol $\wt\forall$ to mean {\em ``for
  almost all.''}
\end{remark}

By Lemma~\ref{l era1} the set $A$ has the
following property
$$\wt\forall (a_0,a_1,x)\in A\times A\times A^c,\ \ Q(a_0,a_2,x)\in A.$$  

We work now in the half space model $\M C\times \M R^+$ of $\M H^3$.
So $S_\infty^2=\M C\cup\{\infty\}$. In that model
$$Q(\infty,a,z)=\alpha(z-a)+a$$ 
where $\alpha=(1+i\sqrt 3)/2$. Again by Fubini's theorem $\wt\forall a_0\in
A,\ \wt\forall (a_1,x)\in A\times A^c$, we have $Q(a_0,a_1,x)\in A$
and we can suppose that this holds for $a_0=\infty$.

In other words,
for almost all $(a,x)\in A\times A^c$ the third vertex of the
equilateral triangle with the first two vertices in $a$ and $x$ is
in $A$. For any $a,x\in\M C$ we call $E_x(a)$ the set of the vertices of
the regular hexagon centered at $x$ and with a vertex in $a$. Then
we have 
\begin{equation}\label{1}
\wt\forall (a,x)\in A\times A^c, \ E_x(a)\subset A
\end{equation}
and in particular 
$\wt\forall (a,x)\in A\times A^c, \ 2x-a\in A.$ Note that $x$ is the
middle-point of the segment between $a$ and $2x-a$.

\begin{lemma}\label{l v2}
For any open set $B\subset\M C$ we have $\mu(A\cap B)>0$.   
\end{lemma}
\proof
Suppose the contrary. Then there exists an open set $B$ 
such that $\mu(A\cap B)=0$. That is,
almost all the points of $B$ are in $A^c$. Moreover, 
from~(\ref{1}) and Fubini's theorem it follows that $\wt\forall x\in
A^c, \wt\forall a\in A, E_x(a)\in A$. 
Therefore there exists a point $x_0\in B$ 
such that a small ball $B_0=B(x_0,r_0)$ is contained in $B$ and 
\begin{equation}\label{2}
\wt\forall a\in A,\ E_{x_0}(a)\in A.
\end{equation}

Since $\mu(A)>0$ then there exists a small ball
$B_1=B(x_1,r_1)$ such that $\mu(A\cap B_1)>0$.
Let $x_2=(x_1+x_0)/2$. If there exists $r>0$ such that $\mu(A\cap
B(x_2,r))=0$, then applying the same argument we can find a point $y$
arbitrarily close to $x_2$ such that~(\ref{2}) holds for $y$. In
particular we get that almost all the points of the
set $C=\{2y-a:\ a\in B_1\cap A\}$ are
in $A$. But if $y$ is close enough to $x_2$ then $C\cap B_0$ has
positive measure, contradicting that $\mu(A\cap B)=0$.

It follows that for all $r_2>0$ we have $\mu(A\cap
B(x_2,r_2))>0$, in particular we choose $r_2<r_0/2$. 
By iterating this construction, we find a sequence of
points $x_n\to x_0$ and radii $r_0/2>r_n>0$ such that 
$\mu(A\cap B(x_n,r_n))>0$. For $n$ large enough this contradicts the
fact that $\mu(A\cap B)=0$.\qed 

\begin{lemma}\label{l v3}
For all $z\in \M C$ we have
\begin{equation}\label{3}
\forall r>0\ \ \mu(B(z,r)\cap A)\geq\frac{1}{2}\mu(B(z,r)).
\end{equation}
\end{lemma}
\proof
From Fubini's theorem, and condition~(\ref{1}), it
follows that for almost all $a\in A$ we have
\begin{equation}\label{4}
\wt\forall x\in A^c, \ E_x(a)\subset A.
\end{equation}
Note that if $(\ref{4})$ holds for $a$, then $(\ref{3})$ holds for $a$.

Let $z\in \M C$. From Lemma~\ref{l v2} it
follows that there exists a sequence $x_n\to z$
such that $(\ref{4})$ (and hence $(\ref{3})$) holds for $x_n$. 
As the function $x\mapsto
\mu(A\cap B(x,r))$ is continuous, then the claim
holds for $z$.\qed

\begin{lemma}
  Let $X\subset \M R^2$ be a measurable set. If there exists
  $\alpha>0$ such that for any ball $B$ 
$$\mu(B\cap X)\geq\alpha\mu(B)$$
then $\mu(\M R^2\setminus X)=0$.
\end{lemma}
This is a standard fact of integration theory
and it follows from Lebesgue's differentiation
theorem (see for example~\cite{Ru}).

From this lemma and Lemma~\ref{l v3} it follows that the set $A$ has
full measure. Since $A=\wb D_\rho^{-1}(p)$ then $\wb D_\rho$ is
constant almost everywhere and Lemma~\ref{l basedue} is proved.\qed

This completes the proof of Theorem~\ref{t ug}.

\section{Corollaries}\label{s cor}

In this section we prove some corollaries that can be useful 
for studying hyperbolic $3$-manifolds.

First we show how
from Theorem~\ref{t ug} one gets a proof of Mostow's rigidity
for non-compact manifolds (see~\cite{P} and~\cite{BCS} for a more
general statement and a different proof).

\begin{teo}[Mostow's rigidity for non-compact manifold]\label{t mstw}
Let $f:M\to N$ be a proper map between two orientable non-compact,
complete hyperbolic $3$-manifolds of finite volume. Suppose that
\vol(M)={\rm deg}(f)\vol(N).
Then $f$ is properly homotopic to a locally isometric covering with the
same degree as $f$.   
\end{teo}
\proof 
Let $\omega$ be the volume form of $N$. For $X=M,N$ let
$\Gamma_X\cong \pi_1(X)$ be the subgroup of ${\rm PSL}(2,\M C)$ such that  
$X=\M H^3/\Gamma_X$.
Let $f_*$ denote both the map induced in homology and the 
representation $f_*:\pi_1(M)\to{\rm PSL}(2,\M C)$. 

First assume that the lift $\wt f:\wt M\to\wt N$ has the cone-property
on the cusps. This implies that $\wt f$ is a pseudo-developing map for $f_*$.
Since $f_*[M]=\deg(f)[N]$ we have
$$
\begin{array}{l}
\vol(M)=\deg(f)\cdot\vol(N)=\langle\omega,\deg(f)[N]\rangle\\
=\langle\omega,f_*[M]\rangle=\langle f^*\omega,[M]\rangle=\vol(f_*).
\end{array}$$

Thus, by Theorem~\ref{t ug} there exists an isometry $\f$ such that 
for any $\gamma\in\Gamma_M$
$$f_*(\gamma)=\f\circ\gamma\circ\f^{-1}.$$
As $\wt M\cong \M H^3$, we consider the 
isometry \f\ as an $f_*$-equivariant map from
$\wt M$ to $\M H^3$. Namely, for any $x\in\M H^3$ and
$\gamma\in\Gamma_M$
$$\f(\gamma(x))=f_*(\gamma)(\f(x)).$$
It follows that $\f$ projects to a locally isometric covering 
$\f:M\to N$ and the convex combination from $\wt f$ to $\f$ projects to
a proper homotopy from $f$ to $\f$. 
{Since the degree of a map is invariant under proper homotopies, then 
$\deg(\f)=\deg(f)$.

We prove now that $f$ is always properly homotopic to a map whose lift
has the cone property on the cusps.
Let $\wt f$ be a lift of $f$. 
For each cusp $N_p=P_p\times[0,\infty)$ let $f_p=\wt
f|_{P_p\times\{0\}}$. Since $f$ is proper it follows that $\wt
f(N_p\times\{\infty\})$ is well-defined. 
Let $F_p:N_p\times[0,\infty)\to\M H^3$ be the map
obtained by coning $f_p$ to $\wt f(N_p\times\{\infty\})$ along geodesic
rays. Let $\wt {f'}$ be the map obtained by replacing, on each cusp $N_p$,
the map $\wt f|_{N_p}$ with the map $F_p$. The map $\wt{f'}$ obviously
has the cone-property on the cusps, and projects to a map 
$f':M\to N$. Moreover, the convex combination from $\wt f$ to
$\wt{f'}$ projects 
to a proper homotopy between $f$ and $f'$.  
\qed

From Theorem~\ref{t ug}, Theorem~\ref{t mstw}, Corollary~\ref{c df}
 and the corresponding statements for compact manifolds, we get the following
statement.
\begin{teo}
  Let M be a complete, oriented hyperbolic $3$-manifold of finite
  volume. Let $\Gamma\cong\pi_1(M)$ be the sub-group of ${\rm PSL}(2,\M C)$
  such that $M=\M H^3/\Gamma$. 
  Let $\rho:\Gamma\to$ {\rm PSL}$(2,\M C)$ be a representation.
  Then $|\vol(\rho)|\leq|\vol(M)|$ and equality holds if and only if
  $\rho$ is discrete and faithful.
\end{teo}

\begin{cor}
Let $M$ be an atoroidal, irreducible, ideally triangulated
$3$-manifold. Let $\z\in\{\M C\setminus\{0,1\}\}^n$ be a solution
of the hyperbolicity equations such that 
$\vol(\z)\neq 0$. Then $M$ is hyperbolic. 
\end{cor}
\proof
This immediately follows from Corollary~\ref{c vg} and 
Thurston's Hyperbolization Theorem~(\cite{Th2}).
\qed 

In~\cite{F3} the notion of {\em geometric} solution of the
hyperbolicity equations is introduced.
Roughly speaking, a geometric solution of the hyperbolicity equations 
for a given ideal triangulation $\tau$ is a choice of moduli 
which is compatible with a global hyperbolic structure on $M$.
In~\cite{F3} it is shown that 
not each solution of hyperbolicity equations is geometric 
(see~\cite{F3} for more details on
algebraic and geometric solutions of hyperbolicity equations).

\begin{cor}\label{c ac}
Let $M$ be a complete hyperbolic $3$-manifold of
finite volume and let $\tau$ be an ideal triangulation of $M$. 
If there exists a solution $\z\in\{\M C\setminus\{0,1\}\}^n$ 
of the hyperbolicity equations for $\tau$, 
then there exists a solution
  $z'$ of hyperbolicity equations that is geometric. Moreover such a
  solution is the one of maximal volume. 
\end{cor}
\proof
Consider a natural straightening of $\tau$, and let $\z'$ be the moduli
induced on $\tau$. By Proposition~\ref{p so}, we have only to prove
that the moduli are not in 
$\{0,1,\infty\}$. Suppose that there is a degenerate tetrahedron
$\Delta_i$. Then
at least two vertices, say $v$ and $w$, of $\Delta_i$  coincide. 

Let $\rho(\z)$ be the holonomy
relative to \z\ and let $D_\z$ be a developing map that is also a
pseudo-developing map for $\rho(\z)$. Then $D_\z$ maps $\Delta_i$ into a
tetrahedron of modulus $z_i$. But by hypothesis, $\z$ is in $\{\M
C\setminus\{0,1\}\}^n$ and so the vertices of $\Delta_i$ are four
distinct points. The last assertion follows from Corollary~\ref{c df}
and Theorem~\ref{t ug}\qed

Corollary~\ref{c ac} tells that, 
once one has a solution $\z\in\{\M C\setminus\{0,1\}\}^n$
of the hyperbolicity equations for a
triangulation $\tau$ of a cusped manifold $M$, 
in order to know if $M$ admits a complete hyperbolic structure of
finite volume, 
it suffices to study the solution of
maximal volume. Namely, if one succeeds to prove that the solution of
maximal volume is geometric, then $M$ is hyperbolic. Conversely, if
one proves that such a solution is not geometric (for example if its
holonomy is not discrete) then $M$ cannot be hyperbolic, and this does
not depend on the chosen triangulation. 

As an example of application
of Corollary~\ref{c ac} we give the following:

\begin{cor}\label{c zeron}
Let $M$ be a cusped $3$-manifold 
equipped with  an ideal triangulation  $\tau$. 
If there exists a solution $\z\in\{\M C\setminus\{0,1\}\}^n$ 
of the hyperbolicity equations for $\tau$, and all the solutions have
zero volume, then $M$ is not hyperbolic. 
\end{cor}

We notice that the hypothesis that all the solutions have
zero volume can be replaced by requiring that the volumes are {\em too small}.
This is because the set of the volumes of the hyperbolic manifolds 
is bounded from below by a positive constant. 

Finally, we obtain another proof of the well-know fact that no
Dehn filling of a Seifert manifold is hyperbolic.
\begin{cor}\label{c fr}
  Let $M$ be a $3$-manifold such that $||(M,\partial M)||=0$ and let $N$
  be a Dehn filling of $M$. Then $N$ is not hyperbolic. 
\end{cor}
\proof Suppose the contrary. Let $\rho$ be the holonomy of the
hyperbolic structure of $N$. From Theorem~\ref{t sv} it follows that
$\vol(\rho)=0$, but from Proposition~\ref{p f2} and 
Corollary~\ref{c nat} it follows that $\vol(\rho)=$ $\vol(N)>0$.\qed

\end{document}